\documentclass[a4paper,11pt]{article}

\usepackage[dvipsnames]{xcolor}
\usepackage[utf8]{inputenc}
\usepackage{comment,enumerate,graphicx}
\usepackage{enumitem,amsmath,amssymb}
\usepackage{skak}
\usepackage{amsthm}
\usepackage{thmtools}
\usepackage{url}
\usepackage{caption}
\usepackage{subcaption}
\usepackage[top=21 mm, bottom=22 mm, left=20 mm, right= 20 mm]{geometry}
\usepackage[hidelinks]{hyperref}
\makeatletter
\newcommand{\labelinthm}[1]{%
   \label{temp#1}
   \protected@write \@auxout {}{\string \newlabel{#1}{{\emph{\ref{temp#1}}}{\thepage}{\emph{\ref{temp#1}}}{temp#1}{}} }%
}
\makeatother
\usepackage{adjustbox}
\usepackage{tikz}
\usetikzlibrary{math,calc,intersections, scopes}
\usetikzlibrary{positioning,arrows,shapes,decorations.markings,decorations.pathreplacing, decorations.pathmorphing,matrix,patterns}
\tikzstyle{vertex}=[circle,draw=black,fill=black,inner sep=0,minimum size=5pt,text=white,font=\footnotesize]
\tikzstyle{redvertex}=[circle,draw=red,fill=red,inner sep=0,minimum size=5pt,text=white,font=\footnotesize]
\definecolor{amber}{rgb}{1.0, 0.75, 0.0}
\definecolor{darkgreen}{rgb}{0.18, 0.7, 0.46}

\declaretheorem[name=Theorem,numberwithin=section]{theorem}

\newtheorem{conjecture}[theorem]{\bf Conjecture}
\newtheorem{question}[theorem]{\bf Question}

\newtheorem*{theorem*}{\bf Theorem}

\theoremstyle{definition}

\newtheorem{definition}[theorem]{\bf Definition}

\newcommand\claimproofend{\renewcommand{\qedsymbol}{$\boxdot$}
\end{proof}
\renewcommand{\qedsymbol}{$\square$}}
\def\eps{\varepsilon}
\def\vx{\textbf{x}}

\newcommand\sizeofpieces{\Large }

\title{Recent progress in graph theory using expansion}
\date{}
\author{Richard Montgomery\thanks{Mathematics Institute, University of Warwick, Coventry, UK. Research supported by the European Research Council (ERC) under the European Union Horizon 2020 research and innovation programme (grant agreement No.\ 947978). Email: {\tt richard.montgomery@warwick.ac.uk}.}}

\begin{document}

\maketitle

\begin{abstract} Graph expansion has long been recognised as an important and desirable property with applications in a wide range of areas in computer science and mathematics. A particular form of expansion known as `sublinear expansion' has recently been used particularly effectively in extremal graph theory, leading to the resolution of many long-standing and notable problems over the last decade and an improved understanding of the structure of sparse graphs. This survey will cover these advances.
\end{abstract}

\section{Introduction.}\label{sec:intro}
Simply put, expanders are graphs with good connectivity properties yet comparatively few edges. Over the last fifty years, expanders have played an important role in a host of areas in computer science and mathematics, from error-correcting codes and efficient communication networks to problems in extremal graph theory. For an overview of many applications of expanders, we recommend the brief introduction by Sarnak~\cite{sarnak2004expander}, the survey of Hoory, Linial and Wigderson~\cite{hoory2006expander} covering many applications in computer science and mathematics, the survey by Lubotzky~\cite{lubotzky2012expander} on further applications in mathematics, and the survey by Krivelevich~\cite{krivelevich2019expanders} on applications in graph theory.
The idea of expanders was first introduced by Bassalygo and Pinsker~\cite{bassalygo1973complexity} in 1973. That expanders exist, e.g.\ that an $n$-vertex graph can have relatively strong connection properties yet only $O(n)$ edges, was originally shown by Pinsker~\cite{pinskerexpands} using essentially a probabilistic argument.
The first explicit examples were given by Margulis~\cite{margulis1973explicit} and by Gabber and Galil~\cite{gabber1979explicit}. This survey will concentrate on recent advances in extremal graph theory using a weak form of expansion known as \emph{sublinear expansion}, which has seen the resolution of a range of long-standing conjectures and other advances. Throughout, we will supress the
more technical details, and happily refer the reader to the recent extensive survey on sublinear expansion by Letzter~\cite{letzter2024sublinear} for more specialist details.

In graph theory, the appearance and use of expanders can roughly be divided into the following four (overlapping) categories:

\noindent\;\;\;\;\textbf{i)} Expanders can be used to provide graphs with desirable properties.

\noindent\;\;\;\;\textbf{ii)} Graphs we are interested in may turn out themselves to be expanders.

\noindent\;\;\;\;\textbf{iii)} An expander can be found in a graph and its properties used to prove results for the original graph.

\noindent\;\;\;\;\textbf{iv)} Expanders may be studied for their own intrinsic interest.

\noindent Examples of \textbf{i)} are that known constructions for expanders can be used to give explicit examples of graphs with large girth and large chromatic number, and that expanders can be used as `templates' in constructive techniques (e.g., in the work in~\cite{bucic2024towards} discussed below). A key example of \textbf{ii)} is that random graphs typically have strong expansion properties (a rich area covered, for example, by the books~\cite{bollobas2011random,janson2011random}), while the `robust' expansion properties of very dense graphs and very dense directed graphs have been used to decompose them into cycles and other graphs (see, for example, the survey by K\"uhn and Osthus~\cite{kuhn2012survey}).
This survey will cover recent results in extremal graph theory in which passing to, and working within, a subgraph with certain expansion properties was central to the proof. Thus, it will fall within category \textbf{iii)}, while also touching upon \textbf{iv)}.

Various definitions of expanders have been used in different settings, and similar properties can be reached through
different routes, e.g.\ by defining an expander from whether random walks on their edges rapidly mix, from their spectral properties, or from combinatorial properties (see~\cite{hoory2006expander} for further details). Here, the sparsity property of our expanders will be indirect. Our expanders will satisfy some combinatorial property which implies a good connection property. We will want to find expanders within sparse graphs. Thus, the classes of expanders we consider will necessarily contain sparse graphs, where the efficacy of our methods will be tested most keenly.

For our continued discussion it will be useful to have some explicit definition of an expander. For this, we will use the following definition of an $\alpha$-expander, where any subset of the vertex set $V(G)$ of a graph $G$ which is not too large has a neighbourhood which is not too small.

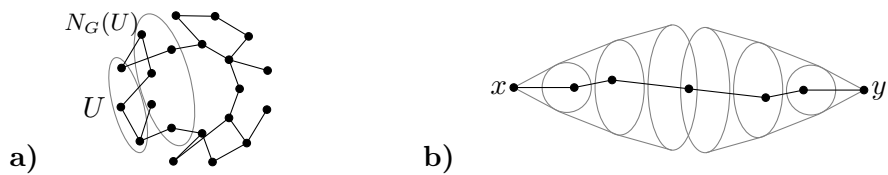
\begin{figure}\centering
\textbf{a)}\;\;\begin{tikzpicture}
\def\radd{1}

\foreach \x in {1,2,3,4,5,6,7,8,9,10,11,12}
{
\coordinate (A\x) at ($(30*\x+14:\radd)$);
\draw [fill] (A\x) circle [radius=0.05cm];
}
\foreach \x in {1,2,3,4,5,6,7,8,9}
{
\coordinate (B\x) at ($(40*\x:0.6*\radd)$);
\draw [fill] (B\x) circle [radius=0.05cm];
}

\begin{scope} [rotate=15]
\draw [black!50]  ($(0.1,0)+(A6)$) circle [x radius=0.2,y radius=0.65];
\end{scope}

\begin{scope} [rotate=15]
\draw [black!50] ($(0.1,0)+(A6)+2*(0.2,0.2)+(0.1,-0.2)+(0.05,0)$) circle [x radius=0.3375,y radius=0.9];
\end{scope}

\draw (A5) -- (A4);
\draw (A5) -- (B3);
\draw (A6) -- (B4);
\draw (A7) -- (B5);
\draw (A7) -- (B6);
\draw (A6) -- (A7);

\foreach \x/\y in {10/8,8/7,8/8,9/7,4/4,3/2,12/1,1/1}
\draw (A\x) -- (B\y);

\foreach \x/\y in {1/2,3/2,11/10,9/10}
\draw (A\x) -- (A\y);

\foreach \x/\y in {1/9,9/8,6/7,2/3,2/1}
\draw (B\x) -- (B\y);



\draw ($(A6)-(0.35,0)$) node {$U$};
\draw ($(A6)+(0.9,0.2)+(-1.15,0.9)$) node {\footnotesize $N_G(U)$};
\end{tikzpicture}\hspace{2cm}\textbf{b)}\;\;
\begin{tikzpicture}
\def\radd{1}

\foreach \x in {1,2,3,4,5,6,7,8,9,10,11,12}
{
\coordinate (A\x) at ($(30*\x+14:\radd)$);
\draw [white] (A\x) circle [radius=0.05cm];
}

\coordinate (x) at ($0.5*(A5)+0.5*(A6)$);
\foreach \n in {1,2,3}
{
\coordinate (x\n) at ($0.5*(A5)+0.5*(A6)+\n*(0.7,0)$);
}

\coordinate (y) at ($0.5*(A12)+0.5*(A11)+(2.7,0)$);
\foreach \n in {1,2,3}
{
\coordinate (y\n) at ($(y)-\n*(0.7,0)$);
}

\draw [black!50] (x) -- ($(x1)+(0,0.33)-(0.1,0.01)$) -- ($(x2)+(0,0.63)-(0.02,0)$) -- ($(x3)+(0,0.83)-(0.01,0)$);
\draw [black!50] (x) -- ($(x1)-(0,0.33)-(0.1,-0.01)$) -- ($(x2)+(0,-0.63)-(0.02,0)$) -- ($(x3)+(0,-0.83)-(0.01,0)$);
\draw [black!50] (y) -- ($(y1)+(0,0.33)-(-0.1,0.01)$) -- ($(y2)+(0,0.63)-(-0.02,0)$) -- ($(y3)+(0,0.83)-(-0.01,0)$);
\draw [black!50] (y) -- ($(y1)-(0,0.33)-(-0.1,-0.01)$) -- ($(y2)+(0,-0.63)-(-0.02,0)$) -- ($(y3)+(0,-0.83)-(-0.01,0)$);

\coordinate (x1pt) at ($(x1)+(0.1,0)$);
\coordinate (x2pt) at ($(x2)+(-0.1,0.1)$);
\coordinate (xypt) at ($0.5*(x3)+0.5*(y3)$);
\coordinate (y1pt) at ($(y1)+(-0.1,0)$);
\coordinate (y2pt) at ($(y2)+(0.1,-0.1)$);

\draw [black!50] (x1) circle [x radius=0.325,y radius=0.33];
\draw [black!50] (x2) circle [x radius=0.325,y radius=0.63];
\draw [black!50] (x3) circle [x radius=0.325,y radius=0.83];
\draw [black!50] (y1) circle [x radius=0.325,y radius=0.33];
\draw [black!50] (y2) circle [x radius=0.325,y radius=0.63];
\draw [black!50] (y3) circle [x radius=0.325,y radius=0.83];

\draw (x) -- (x1pt) -- (x2pt) -- (xypt) -- (y2pt) -- (y1pt) -- (y);

\foreach \n in {x1pt,x2pt,y1pt,y2pt,xypt}
{
\draw [fill] (\n) circle [radius=0.05cm];
}

\draw [fill] (x) circle [radius=0.05cm];
\draw [fill] (y) circle [radius=0.05cm];

\draw ($(x)-(0.2,0)$) node {$x$};
\draw ($(y)+(0.2,0)$) node {$y$};

\end{tikzpicture}
\caption{\textbf{a)} The neighbourhood $N_G(U)$ of a vertex set $U$ in a graph $G$ and \textbf{b)} expanding neighbourhoods iteratively from $x$ and $y$ respectively to find a path from $x$ to $y$.}\label{fig:expansion}
\end{figure}

\begin{definition}\label{defn:linear} An $n$-vertex graph $G$ is an \emph{$\alpha$-expander} if, for each set $U\subset V(G)$ with $|U|\leq n/2$,
\begin{equation}
|N_G(U)|\geq \alpha|U|,\label{eq:expands}
\end{equation}
where $N_G(U)$ is the \emph{neighbourhood of $U$} in $G$, i.e., the set of vertices of $G$ which are not in $U$ and which have at least one neighbouring edge to $U$ (see Figure~\ref{fig:expansion}).
\end{definition}

A basic connection property can be easily seen to follow for any $\alpha$-expander $G$ when $\alpha>0$, by taking iterative neighbourhoods around vertices  (see Figure~\ref{fig:expansion}). For any distinct vertices $x,y\in V(G)$, \eqref{eq:expands} implies that the `balls' $B_G(x)=N_G(x)\cup \{x\}$, $B^2_G(x)=B_G(B_G(x))$, $B^3_G(x)$,{\ldots}  grow at least exponentially in size until they exceed $n/2$ vertices,
as do the balls $B_G(y)$, $B^2_G(y)$, $B_G^3(y)$,\ldots.
Thus, sets in the respective sequences will eventually intersect, giving rise to a path between $x$ and $y$. Moreover, such a path can be taken to have length $O_{\alpha}(\log n)$ (as this holds for arbitrary pairs of vertices $x,y$, the \emph{diameter} of an $\alpha$-expander is $O_{\alpha}(\log n)$).
When $\alpha$ is, say, a large constant, these sets expand in size rapidly, and we can find paths like this with some predictable length, or while avoiding some other vertices in the graph, and this flexibility allows us to start carrying out complex constructions enroute to proving many different results. For more on $\alpha$-expanders, and their properties, we recommend the survey of Krivelevich~\cite{krivelevich2019expanders}. Unfortunately, for each fixed $\alpha>0$, graphs can be fairly dense while still not having any subgraph which is an $\alpha$-expander. We will need to consider what $\alpha$-expanders look like as $\alpha\to 0$.

When $\alpha\to 0$ as $n\to\infty$, this is known as \emph{sublinear expansion}. We will define a specific form of this in Section~\ref{sec:sublinear}, but for now a good question to have in mind is: what can we show in an $n$-vertex $\alpha$-expander if  $\alpha=\log^{-2}n$?
After its introduction by Koml\'os and Szemer\'edi~\cite{komlos1994topological,komlos1996topological} in the mid 1990's, {sublinear expansion} developed slowly in its use as a more general tool. The last decade however has seen a great deal of activity and many advancements. This survey aims to cover many of these advancements.
At the risk of under-serving the many different ideas behind their use of sublinear expansion, we will largely avoid any technicalities and instead only allude to the challenges of working with the weak sublinear expander properties. Interested readers, however, can refer to the recent extensive survey by Letzter~\cite{letzter2024sublinear} which includes excellent summaries of proofs and methods, including for many of the results mentioned here.

In Section~\ref{sec:subdivisions}, we will discuss extremal numbers for subdivisions, the original context in which Koml\'os and Szemer\'edi~\cite{komlos1994topological,komlos1996topological} introduced sublinear expansion, before giving some more details of sublinear expansion in Section~\ref{sec:sublinear}. We will then discuss some problems involving minors (Section~\ref{sec:minors}), cycle lengths (Section~\ref{sec:cyclelengths}), cycle packing (Section~\ref{sec:cyclepacking}), other cycles in sparse graphs (Section~\ref{sec:cycleswithmorproperties}), rainbow cycles (Section~\ref{sec:rainbowcycles}), and cycles in Ramsey theory (Section~\ref{sec:Ramseycycles}). Finally, we will finish by discussing a couple of recent results in which sublinear expansion plays a more unexpected role through its application in an auxiliary graph (Section~\ref{sec:incidental}).


\section{Subdivisions.}\label{sec:subdivisions}
\begin{figure}\centering
\textbf{a)}\;\;\begin{tikzpicture}
\def\radd{1.2}

\foreach \x in {1,2,3,4,5}
{
\coordinate (A\x) at ($(72*\x+18:\radd)$);
\draw [fill] (A\x) circle [radius=0.06cm];
}

\foreach \x/\y in {1/2,2/3,3/4,4/5,1/3,2/4,3/5,2/5,1/4,1/5}
{
\draw (A\x) -- (A\y);
}

\foreach \x/\y/\m/\n in {1/2/3/4,1/2/1/4,1/2/2/4,2/3/1/2,3/4/1/2,4/5/1/2,2/5/1/3,2/5/2/3,2/4/1/6,2/4/2/6,2/4/3/6,2/4/4/6,2/4/5/6}
{
\pgfmathsetmacro\Xmid{1/\n}
\draw [fill] ($\Xmid*\m*(A\x)+(A\y)-\m*\Xmid*(A\y)$) circle [radius=0.06cm];
}
\end{tikzpicture}\hspace{2cm}
\textbf{b)}\;\;\begin{tikzpicture}
\def\radd{1.2}

\foreach \x in {1,2}
\foreach \y in {0,1,2}
{
\coordinate (B\x\y) at ($\x*(1.5,0)+\y*(0,1)$);
\draw [fill] (B\x\y) circle [radius=0.06cm];
}
\foreach \x in {0,1,2}
\foreach \y in {0,1,2}
{
\draw (B1\x) -- (B2\y);
}
\foreach \x/\y/\m/\n in {1/2/3/4,1/2/1/4,1/2/2/4,2/0/0.8/3,2/0/2.2/3,0/0/1/2,2/2/1/3,2/2/2/3}
{
\pgfmathsetmacro\Xmid{1/\n}
\draw [fill] ($\Xmid*\m*(B1\x)+(B2\y)-\m*\Xmid*(B2\y)$) circle [radius=0.06cm];
}
\end{tikzpicture}\hspace{2cm}
\textbf{c)}\;\;\begin{tikzpicture}
\def\radd{1}

\foreach \x in {1,3,4,5,6,8}
{
\coordinate (A\x) at ($(45*\x+22.5:\radd)$);
\draw [fill] (A\x) circle [radius=0.05cm];
}
\foreach \x in {2,7}
{
\coordinate (A\x) at ($(45*\x+22.5:0.9*\radd)$);
\draw [fill] (A\x) circle [radius=0.05cm];
}
\draw (A8) -- (A1) -- (A2) -- (A3) -- (A4) -- (A5) -- (A6) -- (A7) -- (A8);
\draw (A1) -- (A4) -- (A7);
\draw (A2) -- (A4);
\draw (A5) -- (A7)-- (A1);
\draw (A1) to[in=85,out=140] (A3);
\draw (A6) to[in=300,out=-5] (A8);
\end{tikzpicture}
\caption{\textbf{a)} A $K_5$-subdivision, \textbf{b)} a $K_{3,3}$-subdivision, and \textbf{c)} a graph drawn in the plane with no edges crossing, which necessarily contains no $K_5$- or $K_{3,3}$-subdivision.}\label{fig:subdivisions}
\end{figure}
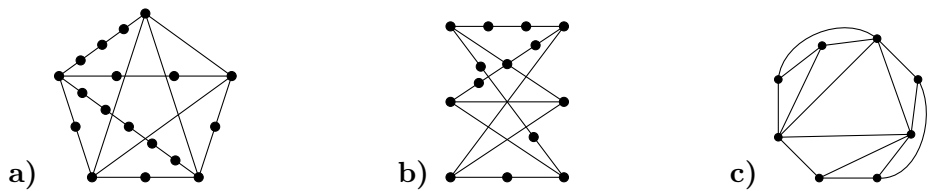

We \emph{subdivide an edge} $xy$ in a graph $G$ by replacing $xy$ with a new vertex whose neighbours are exactly $x$ and $y$. If a graph $G$ can be obtained from another graph $H$ by a sequence of edge subdivisions then we say $G$ is an \emph{$H$-subdivision} (see Figure~\ref{fig:subdivisions}). Subdivisions are sometimes known as \emph{topological minors}, reflecting that many topological properties of graphs are preserved under the subdivision of edges. Indeed, the importance of subdivisions as a concept in graph theory has held since the famous characterisation of planar graphs in 1930 by Kuratowski~\cite{kuratowski1930probleme}. That is, graphs which can be drawn in the plane without edge crossings are exactly those without a subdivision of the complete graph on 5 vertices ($K_5$), or the complete bipartite graph with 3 vertices in each class ($K_{3,3}$) (see Figure~\ref{fig:subdivisions}).

An important class of problems in extremal graph theory asks how many edges a graph on a fixed vertex set must have before the graph necessarily has a particular property of interest. When this property of a graph $G$ is the containment of a copy of another graph $H$ as a subgraph (i.e., edges and vertices can be deleted from $G$ to make it isomorphic to $H$), this is known as a \emph{Tur\'an problem}. The answer has long been known exactly whenever
$H$ is the complete graph $K_n$ (which has $n$ vertices and every possible edge) due to Tur\'an~\cite{turan1941external}, who generalised the earlier result of Mantel~\cite{Mantel} for $n=3$. More generally, when $H$ is not bipartite, the
asymptotic behaviour of the Tur\'an number $\mathrm{ex}(n,H)$, the maximum number of edges an $n$-vertex graph can have yet contain no copy of $H$ is a well-known consequence of the Erd\H{o}s-Stone theorem~\cite{erdos1946structure}. When $H$ is bipartite, $\mathrm{ex}(n,H)$ is known to be $o(n^2)$ (from the Erd\H{o}s-Stone theorem), but determining the correct asymptotic behaviour here is extremely difficult.
For every bipartite graph $H$ containing at least one cycle, it is known there is some $\eps>0$ such that $n^{1+\eps}\leq \mathrm{ex}(n,H)\leq n^{2-\eps}$ (see, for example,~\cite{alon2003turan}). Not only is determining the correct power of $n$ for any such graph $H$ a very difficult problem, finding a wide variety of graphs with different exponents is a fascinating problem. This is the subject of the `rational exponents conjecture' of Erd\H{o}s and Simonovits (see, e.g.,~\cite{erdHos1981combinatorial}), which says that for each rational $\alpha\in (1,2)$ there is some graph $H$ for which $\mathrm{ex}(n,H)=\Theta(n^{\alpha})$. An overview of progress on this problem can be found in the recent work of Conlon and Janzer~\cite{conlon2022rational}, but subdivisions play an important role here as subdividing edges of $H$ should reduce $\mathrm{ex}(n,H)$ if bipartiteness is maintained. Subdivided graphs are good candidates, then, for graphs with different exponents in the asymptotics of their Tur\'an number, as seen for example in the work of Jiang and Qiu~\cite{jiang2020turan}.

If an $n$-vertex graph has $n^{1+o(1)}$ edges then we cannot infer from its number of edges that it contains a copy of any fixed subgraph. In this sparse regime, we may still expect to find some structure in the graph. For example, any $n$-vertex graph with $n$ edges must contain some cycle, and, thus, a subdivision of $K_3$. In such a sparse regime, instead of directly bounding the number of edges we often instead equivalently bound the average degree of a graph, that is, the average number of edges adjacent to the vertices of the graph. Any graph with average degree at least 2 thus contains a subdivision of $K_3$. In 1967, Mader~\cite{mader1967homomorphieeigenschaften} proved that, more generally, for each $t\geq 3$, any graph with sufficiently large average degree (depending only on $t$) must contain a subdivision of $K_t$.
Therefore, we can define $d(t)$ as the smallest $d$ such that every graph with average degree more than $d$ has a $K_t$-subdivision.
Mader's work showed that $d(t)=O(2^{3t})$, but Mader~\cite{mader1967homomorphieeigenschaften}, and Erd\H{o}s and Hajnal~\cite{erdos1969topological}, independently conjectured that $d(t)=O(t^2)$.
In 1994, Bollob\'as and Thomason~\cite{bollobas1996highly} proved this conjecture, very shortly before Koml\'os and Szemer\'edi~\cite{komlos1996topological} improved some of their then-recent work~\cite{komlos1994topological} showing $d(t)=t^{2+o(1)}$ to give a rather different proof, so that we have the following result.

\begin{theorem}\label{thm:extsub} $d(t)=\Theta(t^2)$.
\end{theorem}

Given Theorem~\ref{thm:extsub}, it would be interesting to determine the correct implicit constant.
Due to work by K\"uhn and Osthus~\cite{kuhn2006extremal} in 2006, it is known that $d(t)\leq (\frac{10}{23}+o(1))t^2$. \L uczak  earlier observed that dense random bipartite graphs with a carefully chosen edge density provide an example showing that $d(t)\geq (\frac{9}{64}+o(1))t^2$ (see, for example,~\cite{komlos1996topological}).
Of these bounds, the lower bound seems likely to be nearer the truth.

\begin{conjecture}\label{conj:extsub}
$d(t)= \left(\frac{9}{64}+o(1)\right)t^2$.
\end{conjecture}

The proof of Theorem~\ref{thm:extsub} by Bollob\'as and Thomason~\cite{bollobas1996highly} itself was influential due to its use of, and results on, linked graphs, but even more influential has been Koml\'os and Szemer\'edi's proof due to introduction of what is now known as \emph{sublinear expansion}.


\section{Sublinear expansion.}\label{sec:sublinear}
Having now the prototypical problem for the use of sublinear expansion, we can motivate it further. Suppose we have an $n$-vertex graph $H$ with average degree $d\geq C$, where $C$ is a large constant, and we want to find a $K_4$-subdivision in $H$. Such a subdivision consists of 4 vertices which are pairwise connected by internally-vertex-disjoint paths. Thus, the connectivity properties of expansion make it desirable to find a subgraph of $G$ which is an $\alpha$-expander, for as large a value of $\alpha>0$ as possible.
If all we know is that a graph has large constant average degree, then we cannot take $\alpha$ to be a fixed constant. As shown in the work of~\cite{moshkovitz2018decomposing} mentioned below, the best we could hope for is to take $\alpha=(\log n(\log\log n)^{\Theta(1)})^{-1}$. However, for vertex sets $U$ say of constant size, the bounds on the size of the neighbourhood in \eqref{eq:expands} would then be extraordinarily weak, assuring us only that $U$ has at least one neighbour! We therefore want to have a bound on the size of the neighbourhood of a set that can vary with the size of the set, for which we consider the following definition, due essentially to Koml\'os and Szemer\'edi~\cite{komlos1994topological,komlos1996topological}, in which it would be representative of much of the work we will discuss to think of $\eps$ and $k$ as constants.

\begin{definition}\label{defn:sublinear}
For each $\eps>0$ and $k>0$, an $n$-vertex graph $G$ is an \emph{$(\eps,k)$-expander} if, for each $U\subset V(G)$ with $k\leq |U|\leq n/2$, we have
\begin{equation}\label{eq:expandssublinearlly}
|N_G(U)|\geq \frac{\eps}{\log^2(3|U|/k)}\cdot |U|.
\end{equation}
\end{definition}

Compared to \eqref{eq:expands}, we have replaced $\alpha$ with a function decreasing as $|U|$ increases, so that for large sets our expansion is only by a factor of around $\log^{-2}n$. The key point here is that the expansion in this definition is sufficiently weak that any graph contains an expander whose average degree is at worst only a little smaller than the original graph, as shown by the following result of
Koml\'os and Szemer\'edi~\cite{komlos1994topological,komlos1996topological}.

\begin{theorem}\label{thm:sublinearexists} For each $\delta>0$, there is some $\eps>0$ and $d_0$ such that the following holds for each $d\geq d_0$. Every graph $H$ with average degree $d(H)=d$ contains a subgraph $G$ which is an $(\eps,d/3)$-expander with average degree $d(G)\geq (1-\delta)d(H)$ and which has minimum degree at least $d(G)/2$.
\end{theorem}

The proof of Theorem~\ref{thm:sublinearexists} is by a relatively straight-forward process: starting with $H$, repeatedly, if the graph has too low a minimum degree then we remove a vertex of minimum degree, and if the graph contains a vertex set $U$ which does not expand (i.e., \eqref{eq:expandssublinearlly} fails), then we either \textbf{i)} delete all the vertices in $U$ or \textbf{ii)} delete all the vertices not in $U$ and its neighbourhood. This can be done so that the average degree of the graph does not decrease, unless we have used case \textbf{ii)} but even here the average degree will only decrease slightly (while the number of vertices will go down significantly in this case). The process must end with an expander. To show it has high enough degree still to satisfy Theorem~\ref{thm:sublinearexists}, the decrease in average degree across the process is studied and bounded appropriately. By optimising the various parameters involved it is easy to improve the constant 2 in the power of the logarithm slightly in \eqref{eq:expandssublinearlly}, but it cannot be improved beyond $1+o(1)$, as shown by
Moshkovitz and Shapira~\cite{moshkovitz2018decomposing}. Many variations to Definition~\ref{defn:sublinear} have by now been used, particularly those introducing some `robustness' where the expansion at \eqref{eq:expandssublinearlly} still holds after the arbitrary removal of some small set of edges (including much of the work discussed below). For more on this, we recommend again the survey by Letzter~\cite{letzter2024sublinear}.

Through Theorem~\ref{thm:sublinearexists}, Koml\'os and Szemer\'edi gave us a new paradigm to approach extremal problems in sparse graphs: first we pass to a sublinear expander, then we use its properties to find some desired structure.
The challenge is that such weak expansion properties are very difficult to work with. This can be partially seen by using these conditions to find paths between two vertices $x$ and $y$ in an $(\eps,d)$-expander $G$ with minimum degree at least $d$, similarly to our brief discussion of this in an $\alpha$-expander in Section~\ref{sec:intro}. By a similar argument, we can see there will always be a path between $x$ and $y$ in $G$. However, the weaker expansion rate gives a path of length $O(\log^3n)$ instead of the stronger bound $O_\alpha(\log n)$. An alternative definition of sublinear expansion (see, for example,~\cite{shapira2015small,alon2023essentially}) could improve this to $O(\log n\log\log n)$ but this is essentially best possible (again due to results in~\cite{moshkovitz2018decomposing}). The point we want to make though arises when asking for something just slightly stronger: what if we want a path between $x$ and $y$ with some exact length $\ell$? To make this plausible, say $\ell\approx \log^4n$, so that $\ell$ is larger than any of the diameters of the expanders we have discussed. If $G$ is an $n$-vertex $\alpha$-expander with, say, $\alpha\geq 2$, then the rapid expansion guaranteed by \eqref{eq:expands} in the size of the sets in the sequences $B_G(x)$, $B^2_G(x)$, $B_G^3(x)$,\ldots and $B_G(y)$, $B^2_G(y)$, $B_G^3(y)$,\ldots means that when the sets start to overlap we can take some overlap between the last layers of the respective expansions, so that we know how long the path is that we find between $x$ and $y$. It is not hard to turn this into a precise argument which will find a path of length exactly $\ell$ between $x$ and $y$.
In contrast, while the sets in the sequence $B_G(x)$, $B^2_G(x)$, $B_G^3(x)$,{\ldots} also expand in a sublinear expander, the number of additional vertices at each step may be considerably outnumbered by the previous vertices. I.e., $B_G^i(x)\setminus B_G^{i-1}(x)$ may be much smaller than $B_G^{i-1}(x)$. Therefore, though the sets in the two sequence of increasing balls around $x$ and $y$ will eventually overlap, we will not be able to guarantee the exact length of the path between $x$ and $y$ found by doing so. This drives the difficulty behind finding cycles in sublinear expanders of some exact length, a topic we return to in Section~\ref{sec:cyclelengths}.

Perhaps due to these challenges, the development of tools using sublinear expansion was initially rather slow. After its introduction in the mid 1990's, its application was rather limited in the following years, though notably it was applied in papers by K\"uhn and Osthus~\cite{kuhn2002topological,kuhn2004large} on subdivisions in the early 2000's (and a related notion was used by Matou\v{s}ek~\cite{matouvsek2011number} in discrete geometry). The last decade, however, has seen a great deal of progress. This was initiated at least in large part by a paper of Shapira and Sudakov~\cite{shapira2015small}, first made available in 2012. To discuss this paper, and much of the early subsequent process, we now turn to the study of the extremal function for minors in graphs.


\section{Minors.}\label{sec:minors}
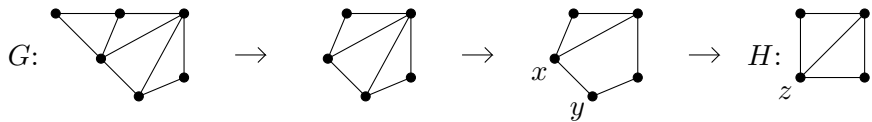
\begin{figure}\centering
\begin{tikzpicture}
\def\radd{0.6}

\foreach \x in {1,2,3,4}
{
\coordinate (A\x) at ($(90*\x+45:\radd)$);
\draw [fill] (A\x) circle [radius=0.06cm];
}
\draw (A1) -- (A2) -- (A3) -- (A4) -- (A1);
\draw (A2) -- (A4);

\foreach \x in {1,2,3,4}
{
\coordinate (B\x) at ($(90*\x+45:\radd)-(3,0)$);
}
\foreach \x in {1,3,4}
\draw [fill] (B\x) circle [radius=0.06cm];

\coordinate (B21) at ($(B2)+(-0.25,0.25)$);
\coordinate (B22) at ($(B2)+(0.25,-0.25)$);
\draw [fill] (B21) circle [radius=0.06cm];
\draw [fill] (B22) circle [radius=0.06cm];
\draw (B3) -- (B4) -- (B1);
\draw (B3) -- (B22) -- (B21) -- (B1);
\draw (B21) -- (B4);

\coordinate (AA) at ($0.5*(A1)+0.5*(B22)+(0.1,0)$);
\coordinate (BB) at ($(AA)-(3,0)$);
\coordinate (CC) at ($(AA)-(6,0)$);

\draw ($(AA)-(9,0)$) node {$G$:};
\draw ($(AA)+(0.8,0)$) node {$H$:};

\foreach \x in {AA,BB,CC}
{
\draw ($(\x)-(0.2,0)$) -- ($(\x)+(0.2,0)$) -- ($(\x)+(0.1,0.1)$);
\draw ($(\x)+(0.2,0)$) -- ($(\x)+(0.1,-0.1)$);
}

\draw ($(B21)-(0.2,0.2)$) node {$x$};
\draw ($(B22)-(0.2,0.2)$) node {$y$};
\draw ($(A2)-(0.2,0.2)$) node {$z$};

\foreach \x in {1,2,3,4}
{
\coordinate (C\x) at ($(90*\x+45:\radd)-(6,0)$);
}
\foreach \x in {1,3,4}
\draw [fill] (C\x) circle [radius=0.06cm];

\coordinate (C21) at ($(C2)+(-0.25,0.25)$);
\coordinate (C22) at ($(C2)+(0.25,-0.25)$);
\draw [fill] (C21) circle [radius=0.06cm];
\draw [fill] (C22) circle [radius=0.06cm];
\draw (C3) -- (C4) -- (C1);
\draw (C3) -- (C22) -- (C21) -- (C1);
\draw (C21) -- (C4);
\draw (C22) -- (C4);

\foreach \x in {1,2,3,4}
{
\coordinate (D\x) at ($(90*\x+45:\radd)-(9,0)$);
}
\foreach \x in {1,3,4}
\draw [fill] (D\x) circle [radius=0.06cm];

\coordinate (D21) at ($(D2)+(-0.25,0.25)$);
\coordinate (D22) at ($(D2)+(0.25,-0.25)$);
\draw [fill] (D21) circle [radius=0.06cm];
\draw [fill] (D22) circle [radius=0.06cm];
\draw (D3) -- (D4) -- (D1);
\draw (D3) -- (D22) -- (D21) -- (D1);
\draw (D21) -- (D4);
\draw (D22) -- (D4);

\coordinate (D5) at ($2*(D1)-(D4)$);
\draw [fill] (D5) circle [radius=0.06cm];
\draw (D4) -- (D5) -- (D21);

\end{tikzpicture}
\caption{$G$ contains $H$ as a minor, for a copy of $H$ can be formed by deleting an edge, then a vertex and then contracting the edge $xy$ into the vertex $z$.}\label{fig:minors}
\end{figure}

For any graph $H$, a graph $G$ \emph{contains an $H$-minor} if a copy of $H$ can be obtained from $G$ by removing edges, removing vertices, and contracting edges (see Figure~\ref{fig:minors}). An early result of Mader~\cite{mader1967homomorphieeigenschaften} shows that an $H$-minor is guaranteed in any graph with sufficiently large average degree (as any $H$-subdivision is an $H$-minor, this is an easier result than the corresponding result of Mader for subdivisions mentioned in Section~\ref{sec:subdivisions}). For each $t\geq 3$, letting $s(t)$ be the least $d$ such that every graph with average degree at least $d$ contains a $K_t$-minor, Mader~\cite{mader1968homomorphiesatze} showed that $s(t)=O(t\log t)$. This is close to the correct behaviour of $s(t)$, for in 1984 Kostochka~\cite{kostochka1984lower} and Thomason~\cite{thomason1984extremal} independently showed that $s(t)=\Theta(t\sqrt{\log t})$. Later, Thomason~\cite{thomason2001extremal} was even able to determine the correct constant $\gamma$ for which $s(t)=(\gamma+o(1))t\sqrt{\log t}$.

Relevant to us is a variation of the extremal problem for minors, introduced by Fiorini, Joret, Theis and Wood~\cite{fiorini2012small}. For each fixed graph $H$, they conjectured that any $n$-vertex graph with average degree only slightly above that required to guarantee an $H$-minor must have moreover an $H$-minor with only $O(\log n)$ vertices.
In the influential work of Shapira and Sudakov~\cite{shapira2015small} mentioned above, it was shown using sublinear expansion that such a minor exists with only $O(\log n\log\log n)$ vertices. Subsequently, the conjecture was proved in full by Montgomery~\cite{montgomery2015logarithmically}, who also showed some analogous results for subdivisions.

Another topic important in the recent development of tools involving sublinear expansion concerns a conjecture of Mader on subdivisions from 1999. Recall from Section~\ref{sec:subdivisions} that some average degree condition of $O(t^2)$ is known to force the existence of a $K_t$-subdivision. Mader~\cite{mader1999extremal} conjectured that an average degree condition only linear in $t$ is sufficient to force the existence of a $K_t$-subdivision in any graph which has no cycles of length 4. (Roughly speaking, this would say that extremal examples for $K_t$-subdivision-free graphs contain some dense structures, which can be avoided with such a $C_4$-free condition.)
Following bounds of K\"uhn and Osthus~\cite{kuhn2004large}, and a related result of Balogh, Liu and Sharifzadeh~\cite{balogh2015subdivisions} for graphs with no cycles of length 6, Mader's conjecture was proved in 2017 by Liu and Montgomery~\cite{liu2017proof}. All of these results use sublinear expansion, where the improvements require many new ideas so that paths can be found while reserving structures to facilitate later constructions (details of which can be found in the papers themselves, as well as to some extent in the survey by Letzter~\cite{letzter2024sublinear}).


\section{Cycle lengths.}\label{sec:cyclelengths}
The discussion of Tur\'an numbers in Section~\ref{sec:subdivisions} implies we cannot expect an $n$-vertex graph with $n^{1+o(1)}$ edges to necessarily have a cycle of any fixed length. However, we can still expect to be able to say something more generally about its cycle lengths. For any graph $G$, the \emph{cycle length set} $\mathcal{C}(G)$ of $G$ (also known as its \emph{cycle spectrum}) is the set of integers $k$ such that $G$ contains a cycle of length $k$.
In 1966, Erd\H{o}s and Hajnal~\cite{erdos1966chromatic} considered the harmonic sum of the cycle lengths of a graph $G$ as a measure of the density of $\mathcal{C}(G)$. In particular, they asked whether the sum $\sum_{\ell\in \mathcal{C}(G)}\ell^{-1}$ can be made arbitrarily large by imposing a bound on the chromatic number $\chi(G)$ of a graph $G$ (i.e., the minimum number of colours required to colour the vertices of $G$ so that every edge is between vertices of different colours). Later, Erd\H{o}s~\cite{erdos1975some} wrote that they felt this should be possible even by imposing a bound only on the average degree of a graph, and this was confirmed in 1984 by Gy\'arf\'as, Koml\'os and Szemer\'edi~\cite{gyarfas1984distribution}, who showed that any graph $G$ with average degree $d$ satisfies $\sum_{\ell\in \mathcal{C}(G)}\ell^{-1}=\Omega(\log d)$.
In work discussed further below, Liu and Montgomery~\cite{liu2023solution} used sublinear expansion to improve this bound to $\sum_{\ell\in \mathcal{C}(G)}\ell^{-1}\geq (1/2-o(1))\log d$.
As anticipated by Erd\H{o}s~\cite{erdos1975some}, the constant $1/2$ is unimprovable here, as if $G$ is the complete bipartite graph with $d$ vertices in each class, then $G$ has average degree $d$ while $\mathcal{C}(G)$ contains only every even number from $4$ to $2d$, so that $\sum_{\ell\in \mathcal{C}(G)}\ell^{-1}=(1/2+o(1))\log d$. Erd\H{o}s~\cite{erdHos1981combinatorial} later conjectured something even stronger: that $\sum_{\ell\in \mathcal{C}(G)}\ell^{-1}$ is minimised over all graphs with average degree at least $d$ by exactly this complete bipartite graph. That this is true when $d$ is sufficiently large is shown in forthcoming work of Milojevi\'c, Montgomery, Pokrovskiy and Sudakov.

It requires many edges to guarantee a cycle of some fixed length. Given instead a sequence $\ell_1,\ell_2,\dots$ of possible lengths, might we need far fewer edges to guarantee a cycle with length in this sequence?
In 1977, Bollob\'as~\cite{bollobas1977cycles} showed that this is possible for infinite arithmetic progressions containing even numbers, confirming a conjecture of Burr and Erd\H{o}s~\cite{erdos1976some}.
We say that a sequence $\ell_1,\ell_2,\ldots$ is \emph{unavoidable with high average degree} if there is some $d$ such that any graph with average degree at least $d$ contains a cycle with length in $\ell_1,\ell_2,\ldots$. In these words, Bollob\'as's result is that any infinite arithmetic progression containing even numbers is unavoidable with high average degree. Note that any sequence which is unavoidable with high average degree must trivially contain infinitely many even numbers, as complete bipartite graphs have very large average degree yet no odd length cycle.

In 2005, Verstra\"ete~\cite{verstraete2005unavoidable} confirmed a conjecture of Erd\H{o}s by showing that there is some sequence of vanishing density which is unavoidable with large average degree. A few years later, Sudakov and Verstra\"ete~\cite{sudakov2008cycle} were able to show that many sparse sequences are unavoidable under the imposition of an average degree condition that grows extremely slowly in the number of vertices. Recently, using sublinear expansion, Liu and Montgomery~\cite{liu2023solution} showed that many sparse sequences are unavoidable with high average degree (i.e., a condition independent of the number of edges in a graph).
In particular, this answered a question of Erd\H{o}s \cite{erdos1984some} by confirming that sufficiently large average degree is enough to guarantee a cycle whose length is a power of 2. On this particular question, Erd\H{o}s and Gy\'arf\'as
(see, for example,~\cite{erdHos1994some}) conjectured that minimum degree 3 should be enough for this, but later thought this may be too strong. Certainly, if the answer to the following question is positive then confirming it would require very substantial new ideas.

\begin{question}\label{qn:mindeg3}
Does every graph with minimum degree at least 3 contain a cycle of length $2^k$ for some $k\geq 2$?
\end{question}

As noted above, in all of these settings we can say nothing meaningful about the odd cycles in a graph. Instead, to guarantee some odd cycles, a natural condition to impose is one on the chromatic number of a graph.
For any graph $G$, let $\mathcal{C}_{\mathrm{odd}}(G)$ be the set of \emph{odd} integers $k$ such that $G$ contains a cycle of length $k$. In 1981, Erd\H{o}s and Hajnal~\cite{erdos1981-20} asked whether $\sum_{\ell\in \mathcal{C}_{\mathrm{odd}}(G)}\ell^{-1}$ can be made arbitrarily large by imposing a bound only on the chromatic number $\chi(G)$ of a graph $G$. In 2011, Sudakov and Verstra\"ete~\cite{sudakov2011cycles} made some progress towards this, by showing that $\sum_{\ell\in \mathcal{C}_{\mathrm{odd}}(G)}\ell^{-1}$ diverges under slightly stronger conditions which are not independent of the number of vertices in $G$. The `odd cycle problem' of Erd\H{o}s and Hajnal was then solved by Liu and Montgomery~\cite{liu2023solution} using their work on cycle lengths in sublinear expanders. They showed that any graph $G$ with chromatic number $k$ satisfies $\sum_{\ell\in \mathcal{C}_{\mathrm{odd}}(G)}\ell^{-1}\geq (1/2-o(1))\log k$, where the constant $1/2$ is unimprovable, as can be seen by considering the complete graph with $k$ vertices. The work in~\cite{liu2023solution} moreover shows that many sparse sequences of odd numbers are unavoidable with high chromatic number (using the corresponding definition to unavoidable sequences with high average degree).

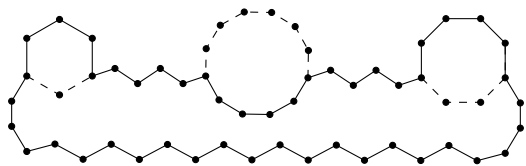
\begin{figure}\centering
\begin{tikzpicture}
\def\radd{0.7}
\def\raddd{0.6}
\def\radddd{0.5}

\foreach \x in {1,2,3,4,5,6,7,8,9,10,11,12}
{
\coordinate (A\x) at ($(30*\x+14:\radd)$);
}

\foreach \x in {1,2,3,4,5,6,7,8}
{
\coordinate (E\x) at ($(3,0)+(45*\x+22.5:\raddd)$);
}
\foreach \x in {1,2,3,4,5,6,7,8}
{
\coordinate (B\x) at ($(3,0)+(45*\x+22.5:\raddd)+(A11)-(E4)+(1.5,0)$);
}

\foreach \x in {1,2,3,4,5,6}
{
\coordinate (D\x) at ($(-3,0)+(60*\x+30:\radddd)$);
}
\foreach \x in {1,2,3,4,5,6}
{
\coordinate (C\x) at ($(-3,0)+(60*\x+30:\radddd)+(A6)-(D5)+(-1.5,0)$);
}

\draw [dashed] (A11) -- (A12) -- (A1) -- (A2) -- (A3) -- (A4) -- (A5) -- (A6);
\draw  (A6) -- (A7) -- (A8) -- (A9) -- (A10) -- (A11);

\draw (B7) -- (B8) -- (B1) -- (B2) -- (B3) -- (B4);
\draw [dashed]  (B4) -- (B5) -- (B6) -- (B7) -- (B8);
\draw (C5) -- (C6) -- (C1) -- (C2) -- (C3);
\draw [dashed]  (C3) -- (C4) -- (C5);

\foreach \x in {1,2,3,4,5,6,7,8,9,10,11,12}
{
\draw [fill] (A\x) circle [radius=0.04cm];
}
\foreach \x in {1,2,3,4,5,6,7,8}
{
\draw [fill] (B\x) circle [radius=0.04cm];
}
\foreach \x in {1,2,3,4,5,6}
{
\draw [fill] (C\x) circle [radius=0.04cm];
}

\def\nn{4}
\foreach \m in {1,...,\nn}
{
\pgfmathsetmacro\Xnew{\nn+1}
\pgfmathsetmacro\Xmid{1/\Xnew}
\pgfmathsetmacro\Xmod{mod(\m,2)}
\coordinate (AA\m) at ($\Xmid*\m*(C5)+(A6)-\Xmid*\m*(A6)+(0,0.1)-\Xmod*(0,0.2)$);
\draw [fill] (AA\m) circle [radius=0.04cm];
}
\foreach \m in {1,...,3}
{
\pgfmathsetmacro\Xnew{\m+1}
\draw (AA\m) -- (AA\Xnew);
}
\draw (C5) -- (AA4);
\draw (AA1) -- (A6);

\def\nn{4}
\foreach \m in {1,...,\nn}
{
\pgfmathsetmacro\Xnew{\nn+1}
\pgfmathsetmacro\Xmid{1/\Xnew}
\pgfmathsetmacro\Xmod{mod(\m,2)}
\coordinate (AAA\m) at ($\Xmid*\m*(A11)+(B4)-\Xmid*\m*(B4)+(0,0.1)-\Xmod*(0,0.2)$);
\draw [fill] (AAA\m) circle [radius=0.04cm];
}
\foreach \m in {1,...,3}
{
\pgfmathsetmacro\Xnew{\m+1}
\draw (AAA\m) -- (AAA\Xnew);
}
\draw (A11) -- (AAA4);
\draw (AAA1) -- (B4);

\coordinate (C3lower) at ($(C3)-(0,1)$);
\coordinate (B7lower) at ($(B7)-(0,1)$);
\draw [fill] (C3lower) circle [radius=0.04cm];
\draw [fill] (B7lower) circle [radius=0.04cm];
\def\nn{16}
\foreach \m in {1,...,\nn}
{
\pgfmathsetmacro\Xnew{\nn+1}
\pgfmathsetmacro\Xmid{1/\Xnew}
\pgfmathsetmacro\Xmod{mod(\m,2)}
\coordinate (AAAA\m) at ($\Xmid*\m*(C3lower)+(B7lower)-\Xmid*\m*(B7lower)+(0,0.1)-\Xmod*(0,0.2)$);
\draw [fill] (AAAA\m) circle [radius=0.04cm];
}
\foreach \m in {1,...,15}
{
\pgfmathsetmacro\Xnew{\m+1}
\draw (AAAA\m) -- (AAAA\Xnew);
}
\draw (C3lower) -- (AAAA16);
\draw (AAAA1) -- (B7lower);

\pgfmathsetmacro\Xmid{1/3}
\coordinate (C31) at ($\Xmid*(C3)+2*\Xmid*(C3lower)+(-0.2,0)$);
\coordinate (C32) at ($2*\Xmid*(C3)+\Xmid*(C3lower)+(-0.2,0)$);
\draw [fill] (C31) circle [radius=0.04cm];
\draw [fill] (C32) circle [radius=0.04cm];
\draw (C3) -- (C32) -- (C31) -- (C3lower);
\coordinate (B71) at ($\Xmid*(B7)+2*\Xmid*(B7lower)+(0.2,0)$);
\coordinate (B72) at ($2*\Xmid*(B7)+\Xmid*(B7lower)+(0.2,0)$);
\draw [fill] (B71) circle [radius=0.04cm];
\draw [fill] (B72) circle [radius=0.04cm];
\draw (B7) -- (B72) -- (B71) -- (B7lower);

\end{tikzpicture}
\caption{A long cycle passing through three short cycles which can be used to adjust the length of the long cycle by changing it to pass instead through the dashed lines in some of the short cycles.}\label{fig:cyclelengthadjuster}
\end{figure}

The underlying mechanism used in~\cite{liu2023solution} to find a cycle of any given long length is very simple. A cycle with approximately the right length is found which moreover contains a structure known as an `adjuster' (see Figure~\ref{fig:cyclelengthadjuster}). This adjuster is simply a collection of `adjusting' cycles attached together into a chain by paths which enter and leave on almost-but-not-quite opposite sides of each adjusting cycle. When choosing a path through these cycles, there are two options for travel through each adjusting cycle. Making appropriate choices allows us to adjust the length of the path through the adjusting cycles, and hence the length of the long cycle.
We mention the structures created in \cite{liu2023solution} in part to emphasise again the challenge behind carrying out constructions within sublinear expanders. In \cite{liu2023solution}, the structures are built very carefully while reserving vertices and other substructures for further connections.

The technical results in~\cite{liu2023solution} do not explicitly guarantee any specific cycle lengths in an $n$-vertex sublinear expander, but sublinear expanders without (roughly) some particular dense substructure are shown to have cycles with every even length in the interval $[\log^7n,n/\log^{12}n]$. The bounds on this interval were not optimised, and the study of which cycles can be guaranteed to exist in a sublinear expander is an interesting subject in its own right.
Sublinear expanders with $n$ vertices may lack cycles with length very close to $n$, for example if they are bipartite graphs with an unbalanced vertex partition. Therefore, it makes sense to consider this question in the setting of regular graphs (where a graph is $d$-regular if every vertex is in exactly $d$ edges). Very recently, Letzter, Methuku and Sudakov~\cite{letzter2025nearly}
considered long cycles in regular sublinear expanders (and applications of this, which can be found in their paper).
Roughly, their work implies here that, when $\eps>0$ is fixed and $n$ is large, for any $d\geq (\log n)^{100}$, any $n$-vertex $d$-regular $(\eps,d)$-expander contains a cycle with length at least $n-n/\log n$. The parameters here may be improvable, and, in particular,
it seems plausible that perhaps even sparser regular sublinear expanders contain a Hamilton cycle (i.e., a cycle through every vertex of the graph), as follows.

\begin{conjecture}\label{conj:cycleinsublinear}
For each $\eps>0$, there exists some $d_0$ such that, for every $d\geq d_0$, every $d$-regular $(\eps,d)$-expander contains a Hamilton cycle.
\end{conjecture}

Even in graphs with a much stronger expansion condition (along the lines of Definition~\ref{defn:linear} with large constance $\alpha$), it is difficult to show that Hamilton cycles exist, as was the subject of an influential conjecture of Krivelevich and Sudakov~\cite{KS:03} from 2003. This conjecture states that, when $d$ is sufficiently large, every $d$-regular expander graph contains a Hamilton cycle, using a definition of expansion via the spectral properties of a  graph. After a long line of research, this was recently confirmed by Dragani{\'c}, Montgomery, Munh{\'a} Correia, Pokrovskiy and Sudakov~\cite{draganic2024hamiltonicity}. Using Definition~\ref{defn:linear}, the more general result of~\cite{draganic2024hamiltonicity} says that there is some constant $C>0$ such that any $n$-vertex $C$-expander, in which moreover any two disjoint vertex sets with size at least $n/C$ have an edge between them, contains a Hamilton cycle.

We conclude this section with a discussion of an interesting problem on cycle length sets where it seems likely the developing techniques involving sublinear expansion will have some import.
The following problem of Erd\H{o}s and Faudree (see, for example,~\cite{erdos1997some}) concerns the number of sets which are realisable as the cycle length set of an $n$-vertex graph. Letting $f(n)$ be the number of subsets $A\subset \{3,\ldots,n\}$ for which there is some $n$-vertex graph $G$ with $\mathcal{C}(G)=A$. What bounds can we achieve on $f(n)$? In addition to the trivial upper bound $f(n)\leq 2^{n-2}$, a construction of Faudree shows that $f(n)\geq 2^{n/2}$. Erd\H{o}s and Faudree asked in particular whether $f(n)=o(2^{n})$ and $f(n)=\omega(2^{n/2})$. In 2004,
Verstra\"ete~\cite{verstraete2004number} showed that the answer to the first question is yes, by showing that $f(n)\leq 2^{n-n^{0.1-o(1)}}$.
As it seems likely that $f(n)=O(2^{n/2})$, answering the second question is potentially extremely difficult, requiring a deep structural result. Very recently, Nenadov~\cite{nenadov2025improved} optimised several components of Verstra\"ete's argument to show that $f(n)\leq 2^{n-n^{0.5-o(1)}}$. The argument requires a structural characterisation of certain classes of graphs, in which the cycle length sets of sparse graphs are particularly relevant. Thus, using sublinear expansion to push these methods further is a promising possibility.


\section{Cycle packing.}\label{sec:cyclepacking}
\begin{figure}\centering
\textbf{a)}\;\;\begin{tikzpicture}
\def\radd{1}

\foreach \x in {1,2,3,4,5,6,7,8,9,10,11,12}
{
\coordinate (A\x) at ($(30*\x+14:\radd)$);
\coordinate (AA\x) at ($(30*\x+14:\radd)$);
\draw [fill] (A\x) circle [radius=0.05cm];
}
\foreach \x in {1,2,3,4,5,6,7,8}
{
\coordinate (B\x) at ($(45*\x+22.5:0.6*\radd)$);
\draw [fill] (B\x) circle [radius=0.05cm];
}
\draw (A1) -- (A2) -- (A3) -- (B2) -- (B1) -- (B8) -- cycle;
\draw (A5) -- (A6) -- (A7) -- (A8) -- (A11) -- (B7) -- (B4) -- (A5);
\draw [densely dashed] (B2) -- (B3) -- (B5) -- (B8) -- (B2);
\draw [densely dotted] (A10) -- (A11) -- (A12) -- (A3) -- (A4) -- (A6) -- (A9) -- (A10);
\draw [densely dash dot] (B6) -- (B7) -- (A10) -- cycle;

\foreach \x in {1,...,12}
{
\coordinate (A\x) at ($(30*\x+14:\radd)-(3.5,0)$);
\draw [fill] (A\x) circle [radius=0.05cm];
}
\foreach \x in {1,2,3,4,5,6,7,8}
{
\coordinate (B\x) at ($(45*\x+22.5:0.6*\radd)-(3.5,0)$);
\draw [fill] (B\x) circle [radius=0.05cm];
}
\draw (A1) -- (A2) -- (A3) -- (B2) -- (B1) -- (B8) -- cycle;
\draw (A5) -- (A6) -- (A7) -- (A8) -- (A11) -- (B7) -- (B4) -- (A5);
\draw  (B2) -- (B3) -- (B5) -- (B8) -- (B2);
\draw  (A10) -- (A11) -- (A12) -- (A3) -- (A4) -- (A6) -- (A9) -- (A10);
\draw  (B6) -- (B7) -- (A10) -- cycle;

\coordinate (AA) at ($0.25*(AA11)+0.25*(AA12)+0.25*(A5)+0.25*(A6)$);
\foreach \x in {AA}
{
\draw ($(\x)-(0.2,0)$) -- ($(\x)+(0.2,0)$) -- ($(\x)+(0.1,0.1)$);
\draw ($(\x)+(0.2,0)$) -- ($(\x)+(0.1,-0.1)$);
}
\end{tikzpicture}\hspace{2cm}
\textbf{b)}\;\;\begin{tikzpicture}
\def\radd{1}

\foreach \x in {1,2,3,4,5,6,7,8,9,10,11,12}
{
\coordinate (A\x) at ($(30*\x+14:\radd)$);
\draw [fill] (A\x) circle [radius=0.05cm];
\coordinate (AA\x) at ($(30*\x+14:\radd)$);
}

\foreach \x in {1,2,3,4,5,6,7,8}
{
\coordinate (B\x) at ($(45*\x+22.5:0.6*\radd)$);
\draw [fill] (B\x) circle [radius=0.05cm];
}

\draw [densely dashed] (A5) -- (A6) -- (A7) -- (A8) -- (A11) -- (B7) -- (B4) -- (A5);
\draw [densely dashed] (B2) -- (B3) -- (B5) -- (B8) -- (B2);
\draw [densely dotted] (A10) -- (A11) -- (A12) -- (A3) -- (A4) -- (A6) -- (A9) -- (A10);
\draw (A1) -- (A12);
\draw (A2) -- (A3);
\draw (B1) -- (B2);
\draw (B6) -- (A9);
\draw (B5) -- (A7);

\foreach \x in {1,...,12}
{
\coordinate (A\x) at ($(30*\x+14:\radd)-(3.5,0)$);
\draw [fill] (A\x) circle [radius=0.05cm];
}
\foreach \x in {1,2,3,4,5,6,7,8}
{
\coordinate (B\x) at ($(45*\x+22.5:0.6*\radd)-(3.5,0)$);
\draw [fill] (B\x) circle [radius=0.05cm];
}
\draw (A5) -- (A6) -- (A7) -- (A8) -- (A11) -- (B7) -- (B4) -- (A5);
\draw (B2) -- (B3) -- (B5) -- (B8) -- (B2);
\draw (A10) -- (A11) -- (A12) -- (A3) -- (A4) -- (A6) -- (A9) -- (A10);
\draw (A1) -- (A12);
\draw (A2) -- (A3);
\draw (B1) -- (B2);
\draw (B6) -- (A9);
\draw (B5) -- (A7);
\coordinate (AA) at ($0.25*(AA11)+0.25*(AA12)+0.25*(A5)+0.25*(A6)$);
\foreach \x in {AA}
{
\draw ($(\x)-(0.2,0)$) -- ($(\x)+(0.2,0)$) -- ($(\x)+(0.1,0.1)$);
\draw ($(\x)+(0.2,0)$) -- ($(\x)+(0.1,-0.1)$);
}
\end{tikzpicture}

\caption{\textbf{a)} An Eulerian graph decomposed into cycles and \textbf{b)} a graph decomposed into cycles and edges.}\label{fig:cycledecomp}
\end{figure}
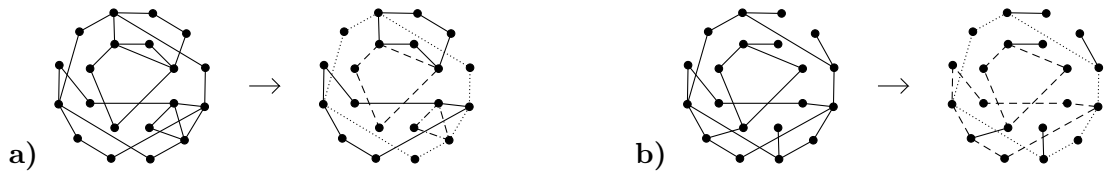

Another area in which sublinear expansion has been effectively used is in packing and decomposition problems, in particular again those problems involving cycles.
It is a simple observation of Veblen~\cite{veblen1912application} that the existence of Euler tours implies that all Eulerian graphs can be decomposed into cycles. That is, if all of the vertex degrees of a graph are even, then the edge set of the graph can be partitioned into cycles (see Figure~\ref{fig:cycledecomp}). A classical result of Walecki from 1892 shows that any Eulerian complete graph can be decomposed into Hamilton cycles (see~\cite{osthus2021extremal}). For much more on the decomposition of dense graphs into different subgraphs, including cycles, we recommend the recent survey of Glock, K\"uhn and Osthus \cite{osthus2021extremal}.

Noting that Walecki's result decomposes Eulerian complete graphs into as few cycles as possible, a natural question is whether every $n$-vertex Eulerian graph has a decomposition into $O(n)$ cycles. The following conjecture of Erd\H{o}s and Gallai~\cite{erdos1966representation} from 1966 is easily seen to be equivalent to this.

\begin{conjecture}\label{conj:EG}
Every $n$-vertex graph has a decomposition into $O(n)$ cycles and edges.
\end{conjecture}

Erd\H{o}s later observed that a construction of Gallai can be improved to show that at least  $(3/2-o(1))n$ cycles and edges may be needed to decompose an $n$ vertex graph. That Conjecture~\ref{conj:EG} remains open is in stark contrast to the corresponding case for decompositions into cycles and \emph{paths}. Here, thanks to an old result of Lov\'asz~\cite{lovasz1968covering}, we have the exactly tight bound that any $n$-vertex graph can be decomposed into at most $\lceil n/2\rceil$ cycles and paths. This follows by an elegant inductive argument, and the lack of a corresponding argument for Conjecture~\ref{conj:EG} indicates it would be a much deeper structural result if true.

When making Conjecture~\ref{conj:EG}, Erd\H{o}s and Gallai~\cite{erdos1966representation} observed that every $n$-vertex graph has a decomposition into $O(n\log n)$ cycles and edges. Such a decomposition can be found by simply iteratively removing a longest cycle from the graph until no cycles (and hence at most $n-1$ edges) remain. A simple undergraduate exercise is that any graph with average degree at least $d\geq 2$ has a cycle with length at least $d/2$. (See also the optimal bounds on this problem in the Erd\H{o}s-Gallai theorem~\cite{gallai1959maximal}.) Tracking, then, the average degree of a graph as longest cycles are removed, it will take the removal of $O(n)$ cycles for the average degree to drop at least by half, at any stage, and thus after the removal of $O(n\log n)$ cycles there will be no more cycles to remove.

Conjecture~\ref{conj:EG} may be considered a stout problem. Despite its inclusion in many of Erd\H{o}s's problem collections, the simple bound of $O(n\log n)$ stood as the state-of-the-art for almost 50 years. Finally, in 2014,  Conlon, Fox and Sudakov~\cite{conlon2014cycle} showed that $O(n\log\log n)$ cycles and edges suffice to decompose any $n$-vertex graph.
A decade later, Buci\'c and Montgomery~\cite{bucic2024towards} reduced this bound further to $O(n\log^*n)$, where $\log^*n$ is the iterated logarithm function, i.e., the least $k$ such that the $k$-fold logarithm of $n$, $\log(\log(\ldots \log (n)))$, is at most $1$.

The improvements in~\cite{conlon2014cycle,bucic2024towards} rely on expansion as a key tool. From Lov\'asz's result quoted above, we know any $n$-vertex graph decomposes into $O(n)$ cycles and paths. The basic approach of Conlon, Fox and Sudakov~\cite{conlon2014cycle} is to set aside vertices independently at random with some small probability, before decomposing the edges not containing any reserved vertex into few cycles and paths. After edge-disjointly joining up the paths into cycles if possible (using the reserved vertices), the argument is then iterated on the edges not covered by cycles so far. If any graph with average degree $d$ can be decomposed into $O(n)$ cycles and a graph of leftover edges of average degree at most (say) $d^{9/10}$, then within $O(\log\log n)$ iterations no cycles will remain, achieving a decomposition into $O(n\log\log n)$ cycles and $O(n)$ edges.
This outline is significantly simplified, and making it work is a difficult task, but essentially Conlon, Fox and Sudakov~\cite{conlon2014cycle} were able to show that it works in graphs satisfying a strong expansion condition. The expansion condition had to be strong enough that some expansion condition involving only the reserved vertices is likely to hold, strong enough that the paths from the decomposition into cycles and paths can be joined up into cycles using the inherited expansion conditions. Notably, the proof in \cite{conlon2014cycle} uses that not only does any graph with enough edges contain an expander, but all except very few of the edges can be decomposed into expanders.

The work of Buci\'c and Montgomery~\cite{bucic2024towards} very roughly followed the outline used by Conlon, Fox and Sudakov, but used sublinear expansion in place of the stronger expansion conditions. This introduces a raft of complications and requires many more ideas, in particular in using sublinear expansion in combination with the random selection of vertices, a difficult proposition made easier by ideas of Tomon~\cite{tomon2022robust} (from work mentioned below).
If Conjecture~\ref{conj:EG} is true, then sublinear expansion may prove useful to prove it, but this appears to be very difficult. Indeed, any reasonable proof that works through a `memoryless' iterative decomposition presumably will at best give a bound of the form $O(n\log^*n)$. Altering cycles found in earlier iterations to reduce the number of cycles used, or decomposing the graph in one step, seems to require a deeper structural understanding of which graphs are difficult to decompose into few cycles.


\section{Cycles with additional properties.}\label{sec:cycleswithmorproperties}
\begin{figure}\centering
\textbf{a)}\;\;\begin{tikzpicture}
\def\radd{1}

\foreach \x in {1,2,3,4,5,6,7,8,9,10,11,12}
{
\coordinate (A\x) at ($(30*\x+14:\radd)$);
\draw [fill] (A\x) circle [radius=0.05cm];
}
\draw [densely dashed] (A5) -- (A9) -- (A2) -- (A11) -- (A6) -- (A3) -- (A5);
\draw (A12) -- (A1) -- (A2) -- (A3) -- (A4) -- (A5) -- (A6) -- (A7) -- (A8) -- (A9) -- (A10) -- (A11) -- (A12);
\end{tikzpicture}\hspace{2cm}
\textbf{b)}\;\;\begin{tikzpicture}
\def\radd{1}

\foreach \x in {1,2,3,4,5,6,7,8,9,10,11,12}
{
\coordinate (A\x) at ($(30*\x+14:\radd)$);
\draw [fill] (A\x) circle [radius=0.05cm];
}
\draw (A12) -- (A1) -- (A2) -- (A3) -- (A4) -- (A5) -- (A6) -- (A7) -- (A8) -- (A9) -- (A10) -- (A11) -- (A12);
\draw [densely dashed] (A3) -- (A5) -- (A8) -- (A11) -- (A1) -- (A3);
\end{tikzpicture}
\hspace{2cm}
\textbf{c)}\;\;\begin{tikzpicture}
\def\radd{1}

\foreach \x in {1,2,3,4,5,6,7,8,9,10}
{
\coordinate (A\x) at ($(36*\x+18:\radd)$);
\draw [fill] (A\x) circle [radius=0.05cm];
}
\draw (A10) -- (A1) -- (A2) -- (A3) -- (A4) -- (A5) -- (A6) -- (A7) -- (A8) -- (A9) -- (A10);
\draw [densely dashed] (A10) -- (A3) -- (A5) -- (A7) -- (A9) -- (A1) -- (A4) -- (A6) -- (A8) -- (A2) -- (A10);
\end{tikzpicture}
\caption{\textbf{a)} Nested cycles, \textbf{b)} nested cycles with no geometric crossing, and \textbf{c)} two cycles nested with each other.}\label{fig:nested}
\end{figure}
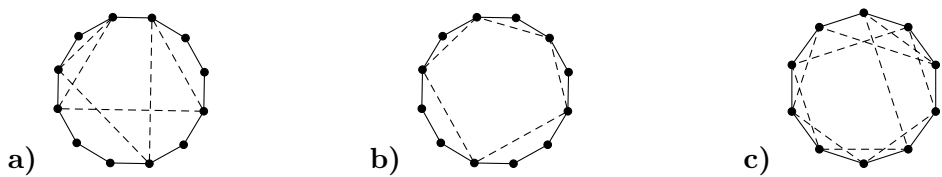

In Section~\ref{sec:cyclelengths}, we discussed what we might be able to say about the length of the cycles in sparse graphs. We now discuss several problems asking whether cycles can be found with additional properties. In all of these problems the extremal number of edges required in an $n$-vertex graph to force the structure considered is $n^{1+o(1)}$, far below the extremal number for any fixed graph.
In 1975, Erd\H{o}s \cite[Problem 29]{erdos1975problems} asked three questions of this kind (for structures depicted in Figure~\ref{fig:nested}). He asked how many edges an $n$-vertex graph can have without, respectively, containing two edge-disjoint cycles which are \textbf{a)} nested, \textbf{b)} nested with no geometric crossings, and \textbf{c)} nested with each other (i.e., with the same vertex set).

The first of these questions was answered fairly promptly by Bollob\'as~\cite{bollobas1978nested}, who showed that any $n$-vertex graph without two edge-disjoint cycles $S_1$ and $S_2$ such that $V(S_1)\subset V(S_2)$ (i.e., such that $S_1$ and $S_2$ are nested) must have $O(n)$ edges. In 1996, 	Chen, Erd\H{o}s and Staton~\cite{chen1996proof} then confirmed a conjecture of Bollob\'as by generalising this result: for each fixed $k$, any $n$-vertex graph without $k$ edge-disjoint cycles $S_1,\dots,S_k$ with $V(S_1)\subset \dots \subset V(S_k)$ has $O(n)$ edges.

In the second question, two edge-disjoint cycles $S_1,S_2$ are \emph{nested with no geometric crossings} if $V(S_1)\subset V(S_2)$ and they can be drawn together in the plane without crossing edges. In 2022,  Fern{\'a}ndez, Kim, Kim and Liu~\cite{gilfernandez2022nested} used sublinear expansion to answer this question, showing that any $n$-vertex graph with no pair of nested cycles with no geometric crossings must have $O(n)$ edges. Interestingly, as the proof takes a shortest cycle in a sublinear expander as the `inner' cycle, the methods seem limited to finding only two such cycles, leaving the following natural conjecture open.

\begin{conjecture}\label{conj:knested} For every $k$, there is some $g(k)$ such that any graph with average degree at least $g(k)$ has $k$ nested cycles with no geometric crossings.
\end{conjecture}

For Erd\H{o}s's third question on nested cycles, in contrast to the other two, it turns out that $\omega(n)$ edges are necessary in an $n$-vertex graph to force two edge-disjoint cycles with the same vertex set. Indeed, due to a surprising construction of Pyber, R\"{o}dl and Szemer\'{e}di~\cite{pyber1995dense} from 1995, there are known to be $n$-vertex graphs with $\Omega(n\log\log n)$ edges yet no 4-regular subgraph (and hence certainly no two edge-disjoint cycles with the same vertex set). For each fixed $k$, the number of  edges required to guarantee a $k$-regular subgraph in an $n$-vertex graph is the topic of the {E}rd{\H{o}}s--{S}auer problem~\cite{erdos1975some} from 1975. This problem was recently resolved in a remarkable work of Janzer and Sudakov~\cite{janzer2023resolution}, who showed that, for each fixed $k$, any $n$-vertex graph without a $k$-regular subgraph has $O(n\log\log n)$ edges. This improved the decades-long best known bound of Pyber~\cite{pyber1985regular} to one that is tight up to the implicit constant.

Similarly strong bounds on the number of edges in an $n$-vertex graph required to guarantee two edge-disjoint cycles with the same vertex set remain elusive, but the first good upper bound was recently shown by Chakraborti, Janzer, Methuku and Montgomery~\cite{chakraborti2025edge}. They showed that, for each fixed $k$, any $n$-vertex graph without $k$ edge-disjoint cycles on the same vertex set must have $n(\log n)^{O(1)}$ edges.
The methods in~\cite{chakraborti2025edge} combine two approaches in extremal graph theory for studying structures in a graph $G$. One, the key theme of this article, is to pass first to a subgraph of $G$ which is a sublinear expander whose average degree is not much smaller than $G$. The second is to pass first to a subgraph of $G$ which is regular (or very nearly regular) without losing too much in terms of the average degree. That this latter approach can be made is closely related to the {E}rd{\H{o}}s--{S}auer problem, where the bound shown by Pyber~\cite{pyber1985regular} is general enough (i.e., effective for large average degree) to find such a regular subgraph without too unreasonable a loss in average degree. Having a regular subgraph then allows the applications of techniques known to work only in regular, or nearly-regular, graphs.

Finding a subgraph which is both regular and a sublinear expander is an interesting challenge. For example, in a $d$-regular graph,  the natural outcome of passing to a sublinear expander as for Theorem~\ref{thm:extsub} will find a subgraph with average degree close to $d$ and degrees between $d$ and slightly less than $d/2$: not sufficiently near regular to be useful for most applications. One key part of the work in~\cite{chakraborti2025edge} is a novel random process to improve the near-regularity of approximately regular graphs, a random process which, if run in a sublinear expander, can be used to show that some weak expansion properties are likely to be retained.
Chakraborti, Janzer, Methuku and Montgomery~\cite{chakraborti2024regular} also showed that this efficient near-regularisation process, along with other ideas, can be used to improve the known results for finding regular subgraphs, in the form of the generalisation of the {E}rd{\H{o}}s--{S}auer problem discussed above. Thanks to the Janzer-Sudakov framework~\cite{janzer2023resolution}, these recent improvements, and a diverse selection of ideas starting with the influential algebraic methods used by Alon, Friedland and Kalai~\cite{alon1984regular} (as used in~\cite{janzer2023resolution}), the minimum number of edges required in an $n$-vertex graph to guarantee a $d$-regular subgraph is now known up to a universal constant multiple for all but a small range of $d$~\cite{chakraborti2024regular}.

\begin{figure}\centering
\textbf{a)}\;\;\begin{tikzpicture}
\def\radd{1}

\foreach \x in {1,2,3,4,5,6,7,8,9,10,11,12}
{
\coordinate (A\x) at ($(30*\x+14:\radd)$);
\draw [fill] (A\x) circle [radius=0.05cm];
}
\draw (A12) -- (A1) -- (A2) -- (A3) -- (A4) -- (A5) -- (A6) -- (A7) -- (A8) -- (A9) -- (A10) -- (A11) -- (A12);
\draw (A1) -- (A9) -- (A3) -- (A6) -- (A2);
\draw (A1) -- (A4) -- (A2) -- (A9) -- (A10);
\draw (A1) -- (A6) -- (A12) -- (A8) -- (A11);
\end{tikzpicture}
\hspace{2cm}
\textbf{b)}\;\;\begin{tikzpicture}
\def\radd{1}

\foreach \x in {1,2,3,4,5,6,7,8,9,10}
{
\coordinate (A\x) at ($(36*\x+18:\radd)$);
\draw [fill] (A\x) circle [radius=0.05cm];
}
\draw (A10) -- (A1) -- (A2) -- (A3) -- (A4) -- (A5) -- (A6) -- (A7) -- (A8) -- (A9) -- (A10);
\draw (A1) -- (A6);
\draw (A2) -- (A7);
\draw (A3) -- (A8);
\draw (A4) -- (A9);
\draw (A5) -- (A10);
\end{tikzpicture}
\hspace{2cm}
\textbf{c)}\;\;\begin{tikzpicture}
\def\radd{1}

\foreach \x in {1,2,3,4,5}
{
\coordinate (A\x) at ($(72*\x+144+90:\radd)$);
\coordinate (B\x) at ($(72*\x+144+90:0.6*\radd)$);
\draw [fill] (A\x) circle [radius=0.05cm];
\draw [fill] (B\x) circle [radius=0.05cm];
\draw (A\x) -- (B\x);
}
\draw (A1) -- (A2) -- (A3) -- (A4) -- (A5) -- (B1);
\draw (B1) -- (B2) -- (B3) -- (B4) -- (B5) -- (A1);
\end{tikzpicture}
\caption{\textbf{a)} A cycle with many chords, \textbf{b)} a cycle with all its diagonals, and \textbf{c)} the graph in \textbf{b)} redrawn as a `twisted cycle of 4-cycles'.}\label{fig:manychords}
\end{figure}
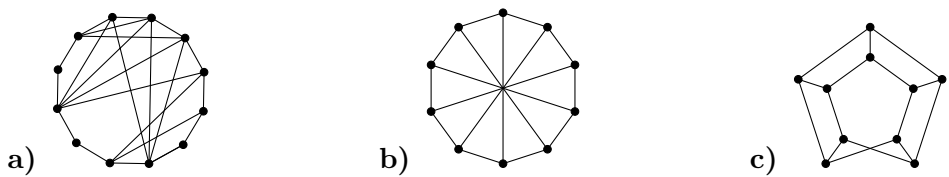

Sublinear expanders with weaker near-regularity properties were also used recently by Dragani{\'c}, Methuku, Munh{\'a} Correia and Sudakov~\cite{draganic2024cycles} for a problem on cycles with many chords. Raised by Chen, Erd\H{o}s and Staton~\cite{chen1996proof} in 1996, this asks how many edges can there be in an $n$-vertex graph without it containing a cycle with at least as many chords as it has vertices (see Figure~\ref{fig:manychords}). By studying random walks in near-regular graphs with weak expansion properties comparable to sublinear expansion, Dragani{\'c}, Methuku, Munh{\'a} Correia and Sudakov~\cite{draganic2024cycles} showed that the answer is $O(n\log^8n)$. This result preceded that on two edge-disjoint cycles with the same vertex set given in~\cite{chakraborti2025edge}, where a weaker bound can be immediately recovered by noting that either of the cycles in such a pair has at least as many chords as edges.

\section{Rainbow Tur\'an problems for cycles and subdivisions.}\label{sec:rainbowcycles}
\begin{figure}\centering
\textbf{a)}\;\;\begin{tikzpicture}
\def\radd{1}

\foreach \x in {1,2,3,4,5}
{
\coordinate (A\x) at ($(72*\x+144+90:\radd)$);
}
\definecolor{colour1}{rgb}{0.93 0.53 0.18};
\definecolor{colour2}{rgb}{0.18 0.55 0.34};
\definecolor{colour3}{rgb}{0.90 0.43 0.16};
\definecolor{colour4}{rgb}{0.0 0.5 1.0};
\definecolor{colour5}{rgb}{0.89 0.0 0.13};

\foreach \x/\y/\col in {1/2/blue,2/3/colour2,3/4/colour3,4/5/colour4,5/1/colour5,3/5/blue,2/4/colour5,1/4/colour2}
{
\draw [ultra thick,\col] (A\x) -- (A\y);
}
\foreach \x in {1,2,3,4,5}
{
\draw [fill] (A\x) circle [radius=0.075cm];
}
\end{tikzpicture}
\hspace{2cm}
\textbf{b)}\;\;\begin{tikzpicture}
\def\radd{1}

\foreach \x in {1,2,3,4,5}
{
\coordinate (A\x) at ($(72*\x+144+90:\radd)$);
}
\definecolor{colour1}{rgb}{0.93 0.53 0.18};
\definecolor{colour2}{rgb}{0.18 0.55 0.34};
\definecolor{colour3}{rgb}{0.90 0.43 0.16};
\definecolor{colour4}{rgb}{0.0 0.5 1.0};
\definecolor{colour5}{rgb}{0.89 0.0 0.13};

\foreach \x/\y/\col in {1/2/blue,2/3/colour2,3/4/colour3,4/5/colour2,5/1/colour5,3/5/blue,2/4/colour5,1/4/colour2}
{
\draw [ultra thick,\col] (A\x) -- (A\y);
}
\foreach \x in {1,2,3,4,5}
{
\draw [fill] (A\x) circle [radius=0.075cm];
}
\end{tikzpicture}
\hspace{2cm}
\textbf{c)}\;\;\begin{tikzpicture}
\def\radd{1}

\foreach \x in {1,2,3,4,5}
{
\coordinate (A\x) at ($(72*\x+144+90:\radd)$);
}
\definecolor{colour1}{rgb}{0.93 0.53 0.18};
\definecolor{colour2}{rgb}{0.18 0.55 0.34};
\definecolor{colour3}{rgb}{0.90 0.43 0.16};
\definecolor{colour4}{rgb}{0.0 0.5 1.0};
\definecolor{colour5}{rgb}{0.89 0.0 0.13};

\foreach \x/\y/\col in {1/2/blue,2/3/colour2,3/4/colour3,4/5/colour4,5/1/colour5}
{
\draw [ultra thick,\col] (A\x) -- (A\y);
}
\foreach \x in {1,2,3,4,5}
{
\draw [fill] (A\x) circle [radius=0.075cm];
}
\end{tikzpicture}
\caption{\textbf{a)} A properly coloured graph, \textbf{b)} a non-properly coloured graph, and \textbf{c)} a rainbow cycle.}\label{fig:colouringconcepts}
\end{figure}
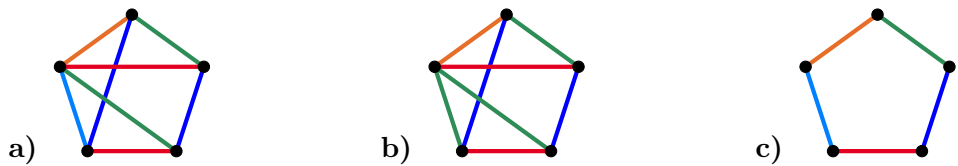

Colouring the edges of graphs is a simple and appealing way to record additional properties or restrictions. The study of edge-coloured graphs, both for their own interest and for their applications, has recently been a particularly active area, as covered by the very recent general survey by Sudakov~\cite{sudakov2024restricted} and the more specialist survey by Pokrovskiy~\cite{Alexeysurvey}. To discuss these problems we recall two concepts (see Figure~\ref{fig:colouringconcepts}). Firstly, that of a \emph{proper colouring}, where each edge of the graph is assigned a colour so that no two edges which contain the same vertex have the same colour. Secondly, we say a colouring of a graph is \emph{rainbow} if each colour appears on at most one edge.
A key example of problems in this area is the so-called `rainbow Tur\'an problem'. Its study was initiated in 2007 by Keevash, Mubayi, Sudakov and Verstra\"ete~\cite{keevash2007rainbow}, who asked, for each fixed graph $H$, how many edges can an $n$-vertex graph have such that it has some proper colouring which has no rainbow copy of $H$. (Note that this is at least $\mathrm{ex}(n,H)$ from Section~\ref{sec:subdivisions}.)
As shown in~\cite{keevash2007rainbow}, the corresponding rainbow Tur\'an number for each $H$ has the same asymptotics as the usual Tur\'an number $\mathrm{ex}(n,H)$ when $H$ is non-bipartite.

As with the classical Tur\'an problem, the problem is most interesting and challenging when $H$ is bipartite. Sometimes the rainbow version of the problem is more approachable. Where our study is seemingly limited by lower bounds on $\mathrm{ex}(n,H)$, we lack graphs without a copy of $H$ with as many edges as we suspect possible. For the rainbow version, good lower bounds can follow from graphs which do have copies of $H$ but have a proper colouring that ensures each of these copies is not rainbow.
A good example is when $H$ is a cycle with $2k$ vertices. Due to Bondy and Simonovits~\cite{bondy1974cycles}, it has been known for more than 50 years that $\mathrm{ex}(n,H)=O(n^{1+1/k})$. Despite generally being expected to tight up to the value of the implicit constant, this has been confirmed only for $k=2,3$ and $5$ (see, for example, the survey by Verstra\"ete~\cite{verstraete2016extremal}). On the other hand, Keevash, Mubayi, Sudakov and Verstra\"ete~\cite{keevash2007rainbow}
were able to show that the corresponding rainbow Tur\'an number for $H$ is $\Theta(n^{1+1/k})$, for all $k\geq 2$.

The most interesting open extremal problem on rainbow cycles is the Tur\'an-style problem for a cycle of any length. Raised as a question by
Keevash, Mubayi, Sudakov and Verstra\"ete~\cite{keevash2007rainbow}, we phrase this as the following conjecture.

\begin{conjecture}\label{conj:rainbowcycle}
Every properly-coloured $n$-vertex graph with no rainbow cycle has $O(n\log n)$ edges.
\end{conjecture}

\noindent This problem has interesting applications to coding theory and to additive combinatorics, for a discussion of which see the survey Sudakov~\cite[Section 2.1]{sudakov2024restricted}.
If true, Conjecture~\ref{conj:rainbowcycle} is correct up to a constant multiple, as seen by the following example of Keevash, Mubayi, Sudakov and Verstra\"ete~\cite{keevash2007rainbow}.
Let $H$ be the graph whose vertex set is the power set of $[k]=\{1,\dots,k\}$, $\mathcal{P}([k])$, where between sets $A$ and $B$ with $A\triangle B=\{i\}$ we put an edge with colour $i$. Thus, $H$ has $n:=2^k$ vertices, and ${kn}/2$ edges, and is properly coloured. Walking around a cycle in $H$ starting with a vertex $A\subset [k]$, we add or remove an element corresponding to the colour of that edge. Thus, as this walk ends back at $A$, the number of edges of each colour in the cycle must be even. In particular, then, $H$ contains no rainbow cycle, but has $n$ vertices and $\Theta(n\log n)$ edges.

From the initial bound towards Conjecture~\ref{conj:rainbowcycle} given by Keevash, Mubayi, Sudakov and Verstra\"ete~\cite{keevash2007rainbow}, there has since been a steady string of improvements. Until recently, the best result was independently shown by Janzer and Sudakov~\cite{janzer2024turan} and~Kim, Lee, Liu and Tran~\cite{kim2022rainbow}, who showed that any properly-coloured graph without an $n$-vertex cycle has $O(n\log^2n)$ edges using an homomorphism counting approach introduced by Janzer~\cite{janzer2023rainbow}. Previously, Das, Lee and Sudakov~\cite{das2013rainbow} had given a bound of the form $n^{1+o(1)}$ using some expansion ideas. With sublinear expansion it is easy to find paths between vertices, but for bounds towards Conjecture~\ref{conj:rainbowcycle} we need methods capable of finding rainbow paths.
This was achieved very effectively in a beautiful paper of Alon, Buci\'c, Sauermann, Zakharov and Zamir~\cite{alon2023essentially}, who showed that every properly-coloured $n$-vertex graph without a rainbow cycle has $O(n\log n\log\log n)$ edges.
The authors of~\cite{alon2023essentially} express some doubts over the veracity of Conjecture~\ref{conj:rainbowcycle}, and it may be that their bound (or something inbetween) is instead the truth. While sublinear expansion may be a crucial tool if further improvements can be made, it is notable that the bound in~\cite{alon2023essentially} matches the best bound on the length of a path that can be found between arbitrary vertex pairs using sublinear expansion (see the relevant discussion in Section~\ref{sec:sublinear}).


A natural generalisation of Conjecture~\ref{conj:rainbowcycle} is to consider the problem for subdivisions of complete graphs more generally, for which the following may be true.

\begin{conjecture}\label{conj:rainbowsubdivision}
 For each $k\geq 3$, every properly-coloured $n$-vertex graph with no rainbow $K_k$-subdivision has $O(n\log n)$ edges.
\end{conjecture}

After an initial bound by Jiang, Methuku and Yepremyan~\cite{jiang2023rainbow}, the first result showing that every properly-coloured $n$-vertex graph with no rainbow $K_k$-subdivision has $n\log^{O(1)}n$ edges was given by Jiang, Letzter, Methuku and Yepremyan~\cite{jiang2021rainbow}. Currently, the best known bound is one of Wang~\cite{wang2022rainbow} of the form $n(\log n)^{2+o(1)}$, proved by optimising a more general argument of Tomon~\cite{tomon2022robust} to this particular setting.
While homomorphism counting approaches give some effective bounds towards Conjecture~\ref{conj:rainbowcycle}, the additional structure required to find rainbow subdivisions more generally seem beyond them and instead all the results towards Conjecture~\ref{conj:rainbowsubdivision} mentioned here follow using expansion in some form. As can be expected, for the state-of-the-art bound of $n(\log n)^{2+o(1)}$ by Wang~\cite{wang2022rainbow}, sublinear expansion is used. While a cursory look at the progress by Alon, Buci\'c, Sauermann, Zakharov and Zamir~\cite{alon2023essentially} towards Conjecture~\ref{conj:rainbowcycle} might suggest a similar bound of $O(n\log n\log\log n)$ can be achieved towards Conjecture~\ref{conj:rainbowsubdivision}, this appears actually to be a very substantial problem. The bound in~\cite{alon2023essentially} follows using a subtle and very clever argument using random selections of the colours which does not seem to be amenable to constructing rainbow subdivisions more generally. This context in particular makes giving a bound of say the form $n\log^{1+o(1)} n$ towards Conjecture~\ref{conj:rainbowsubdivision} a fascinating problem.


\section{Cycles in Ramsey theory.}\label{sec:Ramseycycles}
Another area in which sublinear expansion has been used recently is that of Ramsey theory, and in particular the study of Ramsey numbers involving cycles.
For any graphs $G$ and $H$, the Ramsey number $R(G,H)$ is the least $n$ such that any colouring of the edges of the complete $n$-vertex graph $K_n$ in red and blue must contain a copy of $G$ whose edges are all red or a copy of $H$ whose edges are all blue.
That Ramsey numbers exist in general follows from the foundational result of Ramsey~\cite{RamseyOnAP} in 1930.
An overview of the area can be read in the survey by Conlon, Fox and Sudakov~\cite{conlon2015recent} from 2015, but beyond this the last couple of years have seen truly remarkable progress. In this time extraordinary new bounds have been given (using completely different methods) on the asymptotics of $R(K_t,K_t)$ by Campos, Griffiths, Morris and Sahasrabudhe~\cite{campos2023exponential}, $R(K_4,K_t)$ by Mattheus and Verstra\"ete~\cite{mattheus2024asymptotics} and $R(K_3,K_t)$ by Campos, Jenssen, Michelen and Sahasrabudhe~\cite{campos2025new}.

As it nears the end of its first century, the study of Ramsey numbers of graphs has broadened into a rich and extensive area.  Away from the central goal of improving bounds on Ramsey numbers for complete graphs $R(K_s,K_t)$, when one of the graphs $G$ or $H$ is sparse then we may aim to give better bounds on $R(G,H)$. For example, a classic result of Chv\'atal, R\"odl, Szemer\'edi and Trotter~\cite{chvatal1983ramsey} states that if an $n$-vertex graph $G$ has bounded maximum degree, then $R(G,G)=O(n)$, while much later Lee~\cite{lee2017ramsey} confirmed the Burr-Erd\H{o}s conjecture by showing this holds more widely for graphs $G$ with bounded degeneracy.

Very few Ramsey numbers are known exactly. The value of $R(G,G)$ has long been known when $G$ is a path due to Gerencs\'er and Gy\'arf\'as~\cite{gerencser1967ramsey},
when $G$ is a star (i.e., a graph with a single central vertex which is contained in every edge) due to
Harary~\cite{Harary1972}, and when $G$ is a cycle, due to work by Bondy and Erd\H os~\cite{BONDY197346},  Faudree and Schelp~\cite{faudree1974all}, and Rosta~\cite{rosta1973ramsey}. Much more recently, $R(G,G)$ has been determined when $G$ contains no cycles (i.e., is a tree) and has no very large vertex degrees, by
Montgomery, Pavez-Sign\'e and Yan~\cite{montgomery2025ramsey}.

\begin{figure}\centering
\begin{tikzpicture}[scale=.4]
\def\spacer{3};
\def\Ahgt{1.4};
\def\Bhgt{0.9};


\foreach \m in {1,2,3,4}
{
\coordinate (A\m) at ($\m*(3.5,0)$);
}

\draw[blue,pattern=crosshatch, pattern color=blue!50] ($(A1)+(0,-\Ahgt)$) -- ($(A1)+(0,\Ahgt)$) to[out=30,in=150] ($(A3)+(0,\Ahgt)$) -- ($(A3)-(0,\Ahgt)$) -- cycle;
\draw[blue,pattern=crosshatch, pattern color=blue!50] ($(A2)+(0,-\Ahgt)$) -- ($(A2)+(0,\Ahgt)$) to[out=30,in=150] ($(A4)+(0,\Bhgt)$) -- ($(A4)-(0,\Bhgt)$) -- cycle;
\draw[blue,pattern=crosshatch, pattern color=blue!50] ($(A1)+(0,\Ahgt)$) -- ($(A1)-(0,\Ahgt)$) to[out=-20,in=-160] ($(A4)-(0,\Bhgt)$) -- ($(A4)+(0,\Bhgt)$) -- cycle;

\draw[white,fill=white] ($(A1)+(0,\Ahgt)$) -- ($(A1)-(0,\Ahgt)$) -- ($(A2)-(0,\Ahgt)$) -- ($(A2)+(0,\Ahgt)$) -- cycle;
\draw[white,fill=white] ($(A3)+(0,\Ahgt)$) -- ($(A3)-(0,\Ahgt)$) -- ($(A2)-(0,\Ahgt)$) -- ($(A2)+(0,\Ahgt)$) -- cycle;
\draw[white,fill=white] ($(A3)+(0,\Ahgt)$) -- ($(A3)-(0,\Ahgt)$) -- ($(A4)-(0,\Bhgt)$) -- ($(A4)+(0,\Bhgt)$) -- cycle;

\draw[blue,pattern=crosshatch, pattern color=blue!50] ($(A1)+(0,\Ahgt)$) -- ($(A1)-(0,\Ahgt)$) -- ($(A2)-(0,\Ahgt)$) -- ($(A2)+(0,\Ahgt)$) -- cycle;
\draw[blue,pattern=crosshatch, pattern color=blue!50] ($(A3)+(0,\Ahgt)$) -- ($(A3)-(0,\Ahgt)$) -- ($(A2)-(0,\Ahgt)$) -- ($(A2)+(0,\Ahgt)$) -- cycle;
\draw[blue,pattern=crosshatch, pattern color=blue!50] ($(A3)+(0,\Ahgt)$) -- ($(A3)-(0,\Ahgt)$) -- ($(A4)-(0,\Bhgt)$) -- ($(A4)+(0,\Bhgt)$) -- cycle;

\foreach \m in {1,2,3}
{
\draw [fill=red!60] (A\m) circle [y radius=\Ahgt cm,x radius=1.2cm];
}
\draw [fill=red!60] (A4) circle [y radius=\Bhgt cm,x radius=0.8cm];


\foreach \m in {1,2,3}
{
\draw (A\m) node {$n-1$};
}
\end{tikzpicture}
\caption{Burr's colouring of a complete graph on $(n-1)(\chi(H)-1)+\sigma(H)$ vertices, which contains no red copy of the $n$-vertex connected graph $G$ or blue copy of the graph $H$. When $\chi(H)=4$ this consists of 4 disjoint red complete graphs, 3 with $n-1$ vertices and 1 with $\sigma(H)-1$ vertices, and all the other edges are blue.}\label{fig:burrconstruction}
\end{figure}
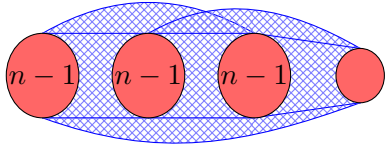

More generally, there are incidences where we could plausibly determine $R(G,H)$ exactly when only one of the graphs $G$ and $H$ is sparse. Due to a construction of Burr~\cite{burr1981ramsey} (see Figure~\ref{fig:burrconstruction}), we know that if $G$ is connected with $n$ vertices then
\begin{equation}
R(G,H)\geq (n-1)(\chi(H)-1)+\sigma(H),\label{eq:Ramseygood}
\end{equation}
where $\sigma(H)$ is the least number of vertices of the same colour over all $\chi(H)$-colourings of $H$. Studying the pairs of graphs for which this bound is tight was initiated in general by Burr and Erd\H{o}s~\cite{burr1983generalizations} in 1983, an area known as `Ramsey goodness'.
Erd\H{o}s had shown that the bound at \eqref{eq:Ramseygood} holds with equality if $G$ is a path and $H$ is a complete graph.
Due to Burr and Erd\H{o}s~\cite{burr1983generalizations} and Erd\H{o}s, Faudree, Rousseau and Schelp~\cite{erdHos1978size}, respectively, it is known that if $G$ is a cycle, or a tree with very low maximum degree, and $H$ has sufficiently few vertices compared to $G$, then the bound at \eqref{eq:Ramseygood} is tight. More recently, Nikiforov and Rousseau~\cite{nikiforov2009ramsey} showed that equality in~\eqref{eq:Ramseygood} holds when $G$ satisfies certain sparseness and separation conditions (see~\cite{nikiforov2009ramsey} for details).

Given these results, attention has turned to qualitative concerns. In 1978, Erd\H{o}s, Faudree, Rousseau and Schelp~\cite{erdHos1978size} made the following conjecture for the Ramsey numbers of cycles versus complete graphs.

\begin{conjecture}\label{conj:Ramseygoodnesscyclecomplete} If $n\geq m\geq 3$, then $R(C_n,K_m)=(n-1)(m-1)+1$.
\end{conjecture}

Nikiforov~\cite{nikiforov2005cycle} showed that this conjecture holds if $n\geq 4m+2$. Perhaps surprisingly, though Conjecture~\ref{conj:Ramseygoodnesscyclecomplete} remains open, it turns out to be true much more broadly when $m$ is large. Indeed,
Keevash, Long and Skokan~\cite{keevash2021cycle} showed that $R(C_n,K_m)=(n-1)(m-1)+1$ holds if $n=\Omega\left(\frac{\log m}{\log\log m}\right)$.
For the Ramsey numbers of cycles $G$ versus graphs $H$ more generally, the following conjecture of Pokrovskiy and Sudakov~\cite{pokrovskiy2020ramsey} (often, but not always, stronger than a related conjecture of Allen, Brightwell and Skokan~\cite{allen2013ramsey}) is appealing.

\begin{conjecture}\label{conj:Ramseygoodnesscycles} There exists $C>0$ such that if $G$ is an $n$-vertex cycle and $n\geq C|H|$ then $R(C_n,H)=(\chi(H)-1)(n-1)+\sigma(H)$.
\end{conjecture}

A corresponding version of Conjecture~\ref{conj:Ramseygoodnesscycles} when $G$ is a path was shown by Pokrovskiy and Sudakov~\cite{pokrovskiy2017ramsey} and when $G$ is a bounded-degree tree by Montgomery, Pavez-Sign\'e and Yan~\cite{montgomery2025ramseybounded}. The best result towards Conjecture~\ref{conj:Ramseygoodnesscycles} is one by Haslegrave, Hyde, Kim and Liu~\cite{haslegrave2023ramsey}, who showed that it holds when $n\geq C|H|\log^4\chi(H)$.
The proof in~\cite{haslegrave2023ramsey} uses sublinear expansion techniques to find a cycle of length $n$ when the graph of red edges satisfies certain expansion properties, using some ideas from the work of Liu and Montgomery~\cite{liu2023solution} discussed above.
Technical barriers from these methods seem to make it hard to make substantial improvements, and a proof of Conjecture~\ref{conj:Ramseygoodnesscycles} would be very interesting, as would be a proof of the remaining unknown cases of Conjecture~\ref{conj:Ramseygoodnesscyclecomplete}.


\section{Expansion in auxiliary graphs.}\label{sec:incidental}
We will finish this survey by discussing briefly two results where sublinear expansion proved very useful when applied not directly to the graph whose properties are under study, but to an auxiliary graph within the proof.
The first of these concerns the extremal number for cycles with all possible diagonals: cycles with even length in which each vertex has a chord connecting it to the vertex furthest from it along the cycle (see Figure~\ref{fig:manychords}).
In 1975, Erd\H{o}s~\cite{erdos1975some} asked for the maximum number of edges an $n$-vertex graph can have without containing a cycle with all its diagonals. Erd\H{o}s~\cite{erdos1975some} observed that there is cycle of length 4 within any cycle with all its diagonals, and that the complete bipartite graph $K_{3,3}$ with 3 vertices in each class is exactly a cycle of length six with all its diagonals. Thus, from known results, the extremal number here is both $\Omega(n^{3/2})$ and $O(n^{5/3})$.

In 2024, Brada{\v{c}}, Methuku and Sudakov~\cite{bradavc2023extremal} showed that the correct answer here is the lower bound, that is, there is some $C>0$ such that any $n$-vertex graph $G$ with at least $Cn^{3/2}$ edges has a cycle with all of its diagonals.
Such a cycle can be thought of as a `cycle of 4-cycles with a twist in it' (see Figure~\ref{fig:manychords}), where a 4-cycle is a cycle with length $4$. To find such a cycle in an $n$-vertex graph $G$ with at least $Cn^{3/2}$ edges, the proof in~\cite{bradavc2023extremal} proceeds by finding a certain cycle in an auxiliary graph whose vertices are edges of $G$ and whose edges correspond to the 4-cycles in $G$. Not only does the proof pass first to an expander subgraph of $G$, an expander containing almost all of the vertices is then found in the corresponding auxiliary graph, using in particular a form of sublinear expansion introduced by Tomon~\cite{tomon2022robust}.

Our second example of sublinear expansion used in an auxiliary graph is in a problem concerning Latin squares.
A \emph{Latin square of order $n$} is an $n$ by $n$ grid filled with $n$ symbols so that each symbol appears exactly once in each row and each column (see Figure~\ref{fig:LS}). A key example is the Cayley table of any group $G$ of order $n$, which forms a Latin square of order $n$. Latin squares have been studied from a mathematical perspective since the work of Euler~\cite{euler}, who considered decompositions of Latin squares into transversals. A \emph{partial transversal} of order $m$ in a Latin square is a set of $m$ cells in the Latin square which share no row, column or symbol. A \emph{transversal} (sometimes called a \emph{full transversal}) is a partial transversal with the same order as the Latin square.

\begin{figure}\centering
\textbf{a)}\;\;\begin{tikzpicture}[define rgb/.code={\definecolor{mycolor}{rgb}{#1}}, rgb color/.style={define rgb={#1},mycolor},scale=0.9]

\def\wi{0.5cm}
\def\n{5}

\draw [white] ($(0,0)$) -- ($(0,-0.5*0.25)$);

\foreach \x in {0,1,...,\n}
\foreach \y in {0,1,...,\n}
{
\coordinate (A\x\y) at ($\n*(0,\wi)+\x*(\wi,0)-\y*(0,\wi)$);
\draw (A\x\y) [thick] rectangle ($(A\x\y)+(\wi,\wi)$);
}

\foreach \x/\y in  {0/0,1/1,2/2,5/4,4/3}
{
\draw (A\x\y) [thick,fill=black!40] rectangle ($(A\x\y)+(\wi,\wi)$);
}

\foreach \x/\y/\symbb in {0/0/0,	0/1/1,	0/2/2,	0/3/3,	0/4/4,	0/5/5,
1/0/1,	1/1/2,	1/2/3,	1/3/4,	1/4/5,	1/5/0,
2/0/2,	2/1/3,	2/2/4,	2/3/5,	2/4/0,	2/5/1,
3/0/3,	3/1/4,	3/2/5,	3/3/0,	3/4/1,	3/5/2,
4/0/4,	4/1/5,	4/2/0,	4/3/1,	4/4/2,	4/5/3,
5/0/5,	5/1/0,	5/2/1,	5/3/2,	5/4/3,	5/5/4
}
{
\draw ($(A\x\y)+0.5*(\wi,\wi)$) node {$\symbb$};
}
\end{tikzpicture}\hspace{2cm}\textbf{b)}\;\;\begin{tikzpicture}[define rgb/.code={\definecolor{mycolor}{rgb}{#1}}, rgb color/.style={define rgb={#1},mycolor},scale=0.9]
\def\wi{0.5cm}
\def\n{5}

\draw [white] ($(0,0)$) -- ($(0,-0.5*0.25)$);

\foreach \x in {0,1,...,\n}
\foreach \y in {0,1,...,\n}
{
\coordinate (A\x\y) at ($\n*(0,\wi)+\x*(\wi,0)-\y*(0,\wi)$);
\draw (A\x\y) [thick] rectangle ($(A\x\y)+(\wi,\wi)$);
}

\foreach \x/\y in  {0/0,4/1,1/2,5/3,2/4,3/5}
{
\draw (A\x\y) [thick,fill=black!40] rectangle ($(A\x\y)+(\wi,\wi)$);
}

\foreach \x/\y/\symbb in {0/0/0,	0/1/1,	0/2/2,	0/3/3,	0/4/4,	0/5/5,
1/0/4,	1/1/0,	1/2/3,	1/3/5,	1/4/2,	1/5/1,
2/0/5,	2/1/3,	2/2/4,	2/3/2,	2/4/1,	2/5/0,
3/0/3,	3/1/4,	3/2/0,	3/3/1,	3/4/5,	3/5/2,
4/0/2,	4/1/5,	4/2/1,	4/3/0,	4/4/3,	4/5/4,
5/0/1,	5/1/2,	5/2/5,	5/3/4,	5/4/0,	5/5/3
}
{
\draw ($(A\x\y)+0.5*(\wi,\wi)$) node {$\symbb$};
}
\end{tikzpicture}
\caption{\textbf{a)} A Latin square of order 6 with a partial transversal of order 5 highlighted and \textbf{b)} a Latin square of order 6 with a (full) transversal highlighted.}\label{fig:LS}
\end{figure}

It is hard to determine whether or not a Latin square has a full transversal, as can be seen in the Hall-Paige conjecture from the 1950's concerning only the very specific case when the Latin square is the Cayley table of a group. The conjecture was eventually proved by a combination of work by Wilcox, Bray and Evans which relied on the classification of finite simple groups (see the survey by Wanless~\cite{wanless2011transversals} for more on this). The natural extremal question for transversals in Latin squares more widely is: how large a partial transversal must there be in any Latin square of order $n$? Towards this, the following combination of conjectures by Ryser~\cite{ryser1967neuere} and Brualdi~(see~\cite{brualdi1991combinatorial}), related to some further conjectures by Stein~\cite{stein1975transversals}, has become the main conjecture in the area, known as the Ryser-Brualdi-Stein conjecture.

\begin{conjecture}\label{conj:RBS}
Every Latin square of order $n$ has a partial transversal with $n-1$ cells and, if $n$ is odd, a full transversal.
\end{conjecture}

After some initial bounds, for almost 40 years the best known bound was one of Shor~\cite{shor} from 1982 (with a proof corrected by Hatami and Shor~\cite{hatamishor} in 2008), showing that every Latin square of order $n$ has a partial transversal with $n-O(\log^2n)$ cells. This was finally improved in 2022 by Keevash, Pokrovskiy, Sudakov and Yepremyan~\cite{KPSY}, who showed that a partial transversal with $n-O(\log n/\log\log n)$ cells always exists in any Latin square of order $n$. The proof in~\cite{KPSY} considers Latin squares in the equivalent form of complete bipartite graphs $K_{n,n}$ which are properly-coloured with $n$ colours, allowing the applications of perspectives and techniques from graph theory (and is a key example of the uses of edge-coloured graphs mentioned in Section~\ref{sec:rainbowcycles}).

Most recently, Montgomery~\cite{montgomery2023proof} showed that, when $n$ is sufficiently large, every Latin square of order $n$ contains a partial transversal with $n-1$ cells. In its lengthy proof,~\cite{montgomery2023proof} uses a variety of tools and novelties, but a critical part of the methods involves decomposing an auxiliary graph into sublinear expanders.
This decomposition is used to find (perhaps degenerate) approximate algebraic structures within a Latin square, and it crucially uses both the good connection properties of sublinear expanders as well as their existence in very sparse graphs. More details can be found in~\cite{montgomery2023proof}, as well as in the survey by the same author~\cite{montgomeryBCCsurvey}, but this is an indication of the unexpected places in which sublinear expansion may prove useful. Sublinear expansion has played a substantial role in the many results covered in this survey, but it is plausible that there is still much progress to be made with this versatile tool.

\section*{Acknowledgments.}
The author would like to thank Matija Buci\'c, Lisa Sauermann, and Benny Sudakov for helpful comments that improved this article.

\bibliographystyle{abbrv}
\bibliography{bibliography}

\end{document}